\documentclass{newclass}

\usepackage{mypackages,Universal,MyAMSEnvironments,graphicx,subfigure,floatflt,citesort}

\newcommand{\minus}{\setminus}

\newcommand{\ra}{\rightarrow}
\newcommand{\R}{\mathbb{R}}

\renewcommand{\epsilon}{\varepsilon}

\renewcommand{\"}[1]{\mathaccent 127 #1}
\renewcommand{\phi}{\varphi}

\newcommand{\cc}{d_{cc}}

\renewcommand{\"}[1]{\mathaccent 127 #1}
\renewcommand{\phi}{\varphi}

\begin{document}

\markboth{Robert K. Hladky and Scott D. Pauls}{Minimal Surfaces in the Roto-translation group}

\title{Minimal surfaces in the roto-translation group with applications to a neuro-biological image completion model}

\author{Robert K. Hladky}

\address{Department of Mathematics, Dartmouth College, 6188 Bradley Hall\\
Hanover, NH 03755 USA\\
robert.k.hladky@dartmouth.edu}

\author{Scott D. Pauls}

\address{Department of Mathematics, Dartmouth College, 6188 Bradley Hall\\
Hanover, NH 03755 USA\\
scott.d.pauls@dartmouth.edu}

\maketitle

\begin{abstract}
We investigate solutions to the minimal surface problem with Dirichlet
boundary conditions in the roto-translation group equipped with a
subRiemannian metric.  By work of G. Citti and A. Sarti, such
solutions are amodal completions of occluded visual data when using a
model of the first layer of the visual cortex.  Using a characterization of  smooth
minimal surfaces as ruled surfaces, we give a method to compute a
minimal spanning surface given fixed boundary data presuming such a
surface exists.  Moreover, we describe a number of obstructions to
existence and uniqueness but also show that under suitable conditions,
smooth minimal spanning surfaces with good properties exist.  Not only
does this provide an explicit realization of the disocclusion process
for the neurobiological model, but it also has application to
contructing disocclusion algorithms in digital image processing.
\end{abstract}

\keywords{Minimal surfaces, Carnot-Carath\'eodory spaces, vision, occlusion}

\ccode{AMS Subject Classification: 53C17,49Q05,92B05}

\section{Introduction}
The study of minimal and isoperimetric surfaces in
Carnot-Carath\'eodory spaces has recently received a good deal of
attention.\cite{CH,CHMY,Cole,DGN,GarNh,GP,LR,LM,Pauls:minimal,Pauls:Obstr,Pauls:cmcgroup,RR}  In
the work cited, the various authors explore the existence, uniqueness
and properties of minimal and isoperimetric problems, finding that, at
least in specific lower dimensional cases, the minimal/isoperimetric
surfaces have a rich geometric structure.  Moreover, recent work of Citti, Manfredini and Sarti\cite{CMS,CS} has provided a link between the method by which the
brain completes missing visual data in the first layer of the visual
cortex (V1) and the solutions to the minimal surface problem in a specific
Carnot-Carath\'eodory space, the roto-translation group, that arises in
a mathematical model of the function of V1.  In this paper, we denote
the roto-translation group by $\mathcal{RT}$.  It is homeomorphic to
$\R^2 \times \sn{1}$ and we will use $(x,y,\theta)$ as coordinates.
Following the construction of Citti and Sarti, we define a
Carnot-Carath\'eodory structure on $\mathcal{RT}$ by distinguishing a
horizontal subbundle, $\mathcal{H}$, given by the span of the following two vector
fields at each point:
\[ X_1 = \cos(\theta) \; \frac{\partial}{\partial x} + \sin(\theta) \;
\frac{\partial}{\partial y}, \; \; X_2=\frac{\partial}{\partial \theta}.\]
These two vector fields bracket generate the tangent bundle and form a
distribution of contact planes in this three dimensional space.
Placing an inner product on $\mathcal{H}$ which makes $\{X_1,X_2\}$ an
orthonormal basis for $\mathcal{H}$, we have the standard
Carnot-Carath\'eodory distance on $\mathcal{RT}$:
\[d_{cc}(a,b) = \inf_{\gamma \in \mathscr{A}}\left \{ \int
<\gamma',\gamma'>^\frac{1}{2} \bigg | \gamma(0)=a, \gamma(1)=b\right \}\]
where $\mathscr{A}$ is the set of absolutely continuous paths whose
derivatives, when they exist, are in $\mathcal{H}$.

In this model, a greyscale image, $I: \Omega \subset \R^2 \ra \R$, has
a representation in $\mathcal{RT}$ given by
\[ \Sigma = \left \{(x,y,\theta) \bigg| \theta = \arctan\left(
    -\frac{I_x}{I_y} \right) \right \}.\]
If a portion of the image is occluded in a domain $\Omega_0 \subset
\Omega$, then Citti and Sarti's model provides a completion of the
occluded region by constructing a minimal spanning surface.  More
precisely, if $c \subset \Sigma$ is the curve in $\mathcal{RT}$
associated to $\partial \Omega_0$, then the completion is given
finding the minimal surface in $\mathcal{RT}$ that spans $c$, i.e. a
minimizer of the perimeter measure.  For a $C^1$ surface, $\Sigma$, given as a
level set of a $C^1$ function $u$, the perimeter is given by
\begin{equation}\label{SA}\mathscr{P}(\Sigma) = \int \sqrt{(X_1 u)^2 +(X_2 u)^2} dA\end{equation}
Moreover we know that such minimal surfaces
satisfy the following partial differential equation:
\[X_1\left( \frac{X_1 u}{\sqrt{(X_1 u)^2 +(X_2 u)^2}} \right ) + X_2
\left ( \frac{X_2 u}{\sqrt{(X_1 u)^2 +(X_2 u)^2}} \right ) =0.\]

In addition to this relationship between the minimal surface problem
and a model for biological image reconstruction, Citti and Sarti
provide a reinterpretation of a number of existing algorithms for
digital inpainting and image completion.  In particular, Citti and
Sarti\cite{CS}, and Citti, Manfredini and Sarti\cite{CMS}, examine the variational models of Ambrosio-Masnou\cite{V:AM}, and a variant of the Mumford-Shah functional
and find that, under suitable interpretation in the roto-translation
group model, minimizing these different functionals is equivalent to
minimizing the standard Carnot-Carath\'edory surface area functional
given in equation \eqref{SA}.  In other words, finding minimizers of these various functionals
is equivalent to solving the Carnot-Carath\'eodory minimal surface problem in
the roto-translation group.  

In light of this unifying theme in the area of vision and image
reconstruction, we explore the minimal
surface problem in a class of groups which include the roto-translation
group, $\mathcal{RT}$.  Citti and Sarti show the divergence form of the minimal  surface equation in $\mathcal{RT}$ and we note that the more general framework of Cheng, Huang,
Malchiodi and Yang\cite{CHMY} shows both the divergence form equation and that the smooth minimal surfaces are ruled surfaces.
We note that the characterization of minimal surfaces as ruled
surfaces has generalizations in Carnot groups with two dimensional
horizontal bundles\cite{Pauls:cmcgroup} and Martinet-type spaces\cite{Cole}.  

In \rfS{PH}, we use the basic form of the
minimal surface equation to explicitly derive the
curves the rule smooth minimal surfaces in the roto-translation group.
Specifically, we show that with respect to 
the Webster-Tanaka connection, $\nabla$, associated to a canonical
pseudo-hermitian structure on $\mathcal{RT}$, the surfaces are
foliated by $\nabla$-geodesics which, for fixed $x_0,y_0,\theta_0,R$ and $\dot{\theta}\neq 0$ take the form:
\begin{equation*}
\begin{split}
x(t) &= x_0+R\sin(\theta(t))\\
y(t) &= y_0 +R\cos(\theta(t))\\
\theta(t) &= \theta_0 + \dot{\theta} t.
\end{split}
\end{equation*}

Thus, we provide a geometric characterization of smooth minimal
surfaces in the roto-translation group which in turn yields an explicit
parametrization for every such minimal surface.  In constrast to the
existing methods of constructing minimal surfaces, which approximate a
minimal surface via a diffusion mechanism, we note that this
parameterization provides a method for constructing exact solutions to
the minimal surface problem.  

Second, we turn to understanding the occlusion problem in
$\mathcal{RT}$.  As demonstrated in the experimental evidence\cite{CS}, in the model of V1 given by the roto-translation group,
representations of image data in $\mathcal{RT}$ potentially contain
different layers of conflicting data due to both modal and amodal
completion of the image.  In light of this finding, we focus on
solving the occlusion problem by finding all possible smooth solutions or
partial solutions of the minimal surface problem with a fixed
boundary.  In \rfS{PP}, we develop a test for determining when
two points on a given curve can be joined by a $\nabla$-geodesic.  We begin with a fixed curve, $c
\subset \mathcal{RT}$, which is the boundary of an occluded region of
the representation of an image in $\mathcal{RT}$ and is parametrized as $c(t)=(\beta(t),\theta(t))$.  For each point, $c(t_0)$ on $c$, we construct the set of other points on $c$
{\em accessible} to $c(t_0)$, denoted $\mathscr{A}(c(t_0),c)$.  To
construct a portion of a smooth minimal surface we simply need construct a
function:
\[ u:D \subset \sn{1} \ra \sn{1}\]
where $D$ is a connected subset of $\sn{1}$ and so that $u(t) \in
\mathscr{A}(c(t),c)$.  We note that to construct a smooth minimal
spanning surface, we must have that $D=\sn{1}$.  For each $t$, this function
give a point $c(u(t))$ connected to $c(t)$ by a rule.  Needless to
say, there are numerous possibilities that occur when attempting to
construct $u$.  In particular, we note that we a guaranteed neither
existence nor uniqueness of such a $u$.  A key tool in the analysis of  $\mathscr{A}(c(t_0),c)$ is the transversality function given by
\[Q(t)=\theta(t)-\phi_\beta(t)\]
where $\phi_\beta$ is defined by the equation
\[\frac{\beta'(t)}{|\beta'(t)|} = (-\sin(\phi_\beta(t)),\cos(\phi_\beta(t)))\]

To further examine this procedure, we simplify the investigation
somewhat and restrict our consideration
to curve $c$ so that the projection of $c$ to $\R^2 \subset \R^2
\times \sn{1} = \mathcal{RT}$ is a circle.  Under this assumption, we are
able to describe a number of different cases ranging from cases where
one can always find such a $u$, cases with multiple $u$ and cases
where no such $u$ exists.  In each of these cases, we give explicit
examples using test image data and give some indication as to the
cause of the various pathologies.  The examples provide a number of
obstruction to the existence and/or uniqueness of smooth minimal
spanning surfaces.

On the positive side, after examining these various cases, we present
a theorem showing sufficient conditions for when a smooth minimal
completion exists.   
\bgT{thm1} Let $I: \R^2 \ra \R$ be an intensity function of an
image with an occlusion given by a circular region $D$.  Further,
suppose $\gamma \in \mathcal{RT}$ is the $\theta$ lift of $\partial
D$ and that the occlusion is completely nondegenerate and occludes
no critical points of $I$.  If $Q'(t) \neq 0$ for $t \in [0,2\pi]$
then there exists a minimal spanning surface of $\gamma$ where the projection of each rule of the surface lies in the interior of $D$.  Moreover, if $Q'(t) < 0$ for $t \in [0,2\pi]$ then the projection of this spanning surface to the xy-plane is surjective onto the occluded region.
\enT

An occlusion is completely nondegenerate if there are no critical
points of $I$ on $\partial D$, only a finite number of critical
points in the interior of $D$ and an angle function can be extended
continuously across those critical points (see below for a more
precise definition).

We again emphasize that, in such a case, the construction of such a surface is significantly less
computationally intensive than the iterative approximative method
used by Citti and Sarti\cite{CS}.  We expect that similar gains can be achieved
with respect to the other models mentioned above such as the
Ambrosio-Masnou and elastica methods.  In this direction, we note that
the authors\cite{HP3} use a discrete version of the method used in this theorem to provide a new algorithm for disocclusion in the context of digital
image reconstruction.   Moreover, we expect that this
method will have application to neurobiology:  by explicitly
constructing completions of images, we will be able to provide
testable hypotheses for neurobiological function of V1.

\section{Modeling V1 via the roto-translation group}\label{model}
In this section, we review the basic biological findings describing
the function of V1 and describes a mathematical model of V1.
Moreover, we describe the connection, provided by Citti and Sarti\cite{CS}, between minimal surfaces in the model space and solutions
to the problem of amodally completing regions of occluded image data.  

\begin{figure}
\begin{center}
\includegraphics[height=3in,width=3in]{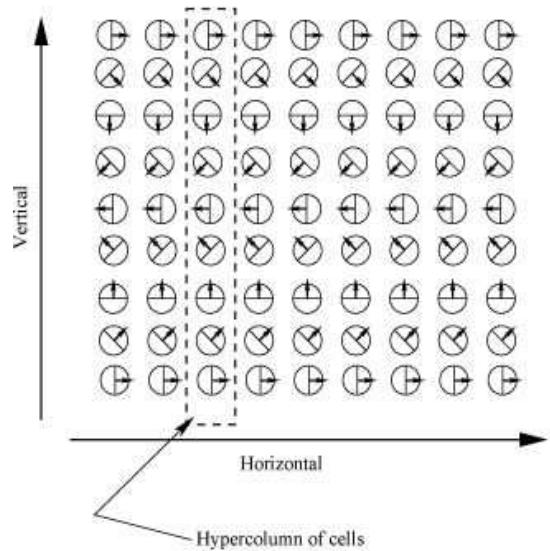}
\caption{A schematic of the structure of V1}\label{rtcols}
\end{center}
\end{figure}

Over the past several decades, the function and operation of the first
layer of the visual cortex, V1, has become increasingly clear.  Early
research showed that V1 contains so-called simple cells that are
sensitive to, among other things, brightness gradients with a
particular orientation.  These cells are arranged in columns sharing
the same orientation preference\cite{V:HW62,V:HW77} and the columns are arranged in hypercolumns which represent all possible orientations.  This
view was further explored and modeled mathematically,\cite{V:Hoffman,V:PT}
where the authors modeled the hypercolumnar cell structure using a contact manifold.  The contact model is based on a simplifying assumption that treats each column as a point, ignoring the column structure to focus on the hypercolumn structure.  Mathematically, they use the manifold $\R^2
\times \sn{1}$ to model the hypercolumn structure by placing a circle of directions above each point $(x,y)
\in \R^2$.  Each point $(x,y,\theta)$ represents a column of cells associated to an $(x,y)$ point of retinal data, all of which are
attuned to the orientation give by the angle $\theta$.  See figure
\ref{rtcols} for a schematic of the hypercolumnar structure.

Early assumptions that cortical connectivity should run mostly
vertically along
the hypercolumns and be severely restricted in horizontal directions, while
supported by some research, was contradicted by later evidence which
showed that there is ``long range horizontal'' connectivity in the cortex.  These experiments (see for example, Gilbert et al\cite{V:GDIW})  indicated that
horizontal connections are made between cells in different hypercolumns of
similar orientation preference.  Moreover, experimental evidence
showed that there is a stronger preference for communication between
cells of not only similar orientation preference but for ones that lie
(roughly) along the axis corresponding to the shared orientation.  In
other words, using the notation of the above model, if
$(x,\theta)$ and $(y,\theta)$ are points in different hypercolumns
with the same angle preference $\theta$, communication between the
cells is preferred if the direction $\theta$ corresponds with the
direction of the vector from $x$ to $y$ in $\R^2$.  This evidence points towards a
geometric structure in this layer where communication between adjacent
cells is allowable in certain directions, vertically and between cells
in different hypercolumns of similar orientation sensitivity, and
vastly restricted in all other directions.  This type of situation has
been studied in a variety of settings including, for example, control
theoretic problems where the degrees of freedom at a particular point
are restricted.  

Petitot and Tondut\cite{V:PT,V:P} incorporate the
these biological findings into their model by introducing a contact
structure on $\R^2 \times \sn{1}$ via the one form $\omega=dx-\theta dy$
and introduce a sub-Riemannian metric associated to the contact
two-plane distribution to encode the geometry of the model of V1.  The
plane field given by the kernel of $\omega$, $span \{ \partial_\theta,
\partial_y + \theta \partial_x\}$, corresponds to the space
of allowable directions at each point.  Notice that these are
precisely the vertical direction and the direction which links cells
in different hypercolumns with the same $\theta$ value.  Citti and Sarti\cite{CS} use the following explicit realization
of the roto-translation group, $\mathcal{RT}$:
\begin{itemize}
\item $\mathcal{RT}$ is diffeomorphic to $\R^2 \times \sn{1}$ with
  coordinates $(x,y,\theta)$.  
\item The following three vector fields span the tangent space at each
  point:
\begin{equation}\label{rtrep}
\begin{split}
 X_1 &= \cos(\theta) \; \frac{\partial}{\partial x} + \sin(\theta) \;
  \frac{\partial}{\partial y}\\
X_2 &= \frac{\partial}{\partial \theta}\\
X_3 &= -\sin(\theta) \; \frac{\partial}{\partial x} + \cos(\theta) \;
  \frac{\partial}{\partial y}
\end{split}
\end{equation}
We note that $[X_2,X_1]=X_3$ and $[X_2,X_3]=-X_1$.
\item For an image $I:D \subset \R^2 \ra \R$, its representation, $\Sigma(I)$, in
  $\mathcal{RT}$ is given by $(x,y,\theta(x,y))$ where $(x,y) \in D$
  and $\theta$ is given by
\[ \frac{\nabla I}{|\nabla I|} = (-\sin(\theta),\cos(\theta))\]
\end{itemize}
We note that this is an explicit realization of the
model described above and matches with the biological evidence concerning horizontal
connectivity.  The contact subbundle in this presentation is simply
$span\{X_1,X_2\}$.  It is a direct calculation that this subbundle
gives a contact structure and, placing an inner product on this
subbundle making $\{X_1,X_2\}$ orthonormal, we have a standard
sub-Riemannian metric on $\mathcal{RT}$ (see the next section for a
precise definition).  One of the main contributions of Citti and
Sarti's adaptation of the cortical model is the use of an explicit
lifting function that transforms retinal data into a surface in the
cortex and allows direct use of the Carnot-Carath\'eodory structure.  

\begin{figure}
\begin{center}
\includegraphics[height=4in,width=2.75in]{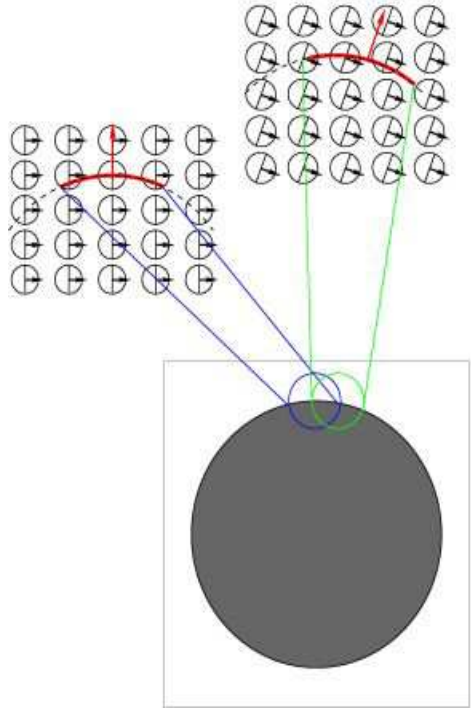} 
\caption{Example of lifting an image to $\mathcal{RT}$}\label{rtpic1}
\end{center}
\end{figure}

To explore this further, we investigate the application
of the model to an image.  Using the representation above, we see that the
direction given by the
angle $(\cos(\theta),\sin(\theta))$ points in a direction
perpendicular to the gradient of $I$, i.e. a direction tangent to the
level sets of $I$.  Thus, if nearby points have the same intensity and
thereby lie on the same level sets of $I$, their $\theta$ representations will
be the same.  As $\theta$ denotes a position in a hypercolumn of
cells over the point $(x,y)$, this echoes the biological finding
horizontal communication occurs between cells of similar orientation
specificity and the property that the representation respects level lines reflects the biological principle that communication between point $(x_0,y_0,\theta)$ and $(x_1,y_1,\theta)$ is permitted if $\theta$ points in the same direction as the vector from $(x_0,y_0)$ to $(x_1,y_1)$.  In figure \ref{rtpic1}, we give a schematic of lifting a
simple image to $\mathcal{RT}$.  Represented are two layers of cells
of similar orientation preference but different hypercolumns and two
sections of the image that could plausibly be lifted to those layers.

Citti, Sarti and Manfredini\cite{CMS} provide a link between mean
curvature flow in Riemannian approximates of the sub-Riemannian space
and the celebrated Mumford-Shah functional.  Further, Citti and Sarti\cite{CS} provided a link between the flow mechanism and several
other models including elastica methods\cite{V:NMS} and the Ambrosio-Masnou model\cite{V:AM}.  Moreover, they showed that surfaces that were
asymptotically stationary under the flows in the Riemannian
approximates to $\mathcal{RT}$ are minimal surfaces in the
sub-Riemannian roto-translation group.  As discussed in the
introduction, minimal surfaces in
Carnot-Carath\'eodory spaces have been examined in some generality\cite{DGN,GarNh,Pauls:cmcgroup} as well as in more
restricted settings.\cite{CH,CHMY,GP,Pauls:Obstr,Pauls:minimal}

Citti and Sarti\cite{CS} use the model of V1 to investigate  the problem of amodally filling in an image when a portion of the image is missing due to
occlusion or some other factor.  To use the model, Citti and Sarti
took image data with a portion deleted and used an approximation of the
flow described above to find a minimal spanning surface that ``fills''
the hole in the image.  In other words, they find solutions to the
minimal surface problem subject to Dirichlet boundary conditions.  In
the next sections, we will examine the minimal surface problem in the
setting of a class of sub-Riemannian spaces which include the
roto-translation group.  The main goal of this paper is to provide a description of such minimal surfaces and to describe method by which they can be constructed in particularly a simple manner.

\section{Notation}\label{notation}
In this section, we fix the basic notation used throughout the paper.
Let $\mathcal{G}$ be a topologically three dimensional one step graded
Lie group.  In other words, the Lie algebra of left invariant vector
fields $\mathcal{V}$ splits as 
\[ \mathcal{V} = \mathcal{V}_0 \oplus \mathcal{V}_1, \; \; \text{dim}
\mathcal{V}_0 = 2, \; \mathcal{V}_1 = [ \mathcal{V}_0,\mathcal{V}_0]\]
Moreover, we assume the following
\begin{itemize}
\item  $\mathcal{G}$ is equipped with a Riemannian
metric $g$, which we sometimes denote in inner product notation by $\langle
\cdot, \cdot\rangle$ and which makes the grading orthogonal.  
\item $\nabla^{LC}$ is the Levi-Civita connection associated to $g$

\item {\em Assumption 1:} $\mathcal{G}$ is said to satisfy assumption
  1 if 
\[[\mathcal{V}_0,\mathcal{V}_1] \subset \mathcal{V}_0\]
\item We define a Carnot-Carath\'eodory distance on $\mathcal{G}$ by
\[ \cc(x,y) = \inf_{\gamma \in \mathcal{A}} \left \{ \int
  <\gamma',\gamma'>^\frac{1}{2} \bigg | \gamma(0)=x,\gamma(1)=y\right
  \}\]where $\mathcal{A}$ is the set of all absolutely continuous
  paths that, where their derivatives are defined, have $\gamma' \in
  \mathcal{V}$.  
\end{itemize}
We note that if $\mathcal{G}$ is nilpotent, then $G$ is known as a {\em
    Carnot group}.  We review two special examples.
\begin{exm}  {\em The Heisenberg group}.  The topologically three dimensional Heisenberg group, $\mathbb{H}$,  is one of the simplest nonabelian nilpotent Lie groups.  As a smooth manifold, it is diffeomorphic to $\R^3$.  Using the terminology above, we have the Lie algebra given as
\[\mathfrak{h} = \mathcal{V}_0 \oplus \mathcal{V}_1\]
with $\mathcal{V}_0 =span \{X_1,X_2\}$ and $\mathcal{V}_1 = span \{X_3\}$ where there is a single nontrivial bracket operation, $[X_1,X_2]=X_3$.  We note that, as $\mathbb{H}$ is nilpotent, it is a Carnot group.
\end{exm}
\begin{exm} {\em The Roto-Translation group}.  As the roto-translation group, $\mathcal{RT}$, appears centrally in the model of visual processing in V1, we review its abstract structure.  As a smooth manifold, it is diffeomorphic to $\R^2 \times \sn{1}$ and, using the terminology above, we have the Lie algebra given as
\[\mathfrak{h} = \mathcal{V}_0 \oplus \mathcal{V}_1\]
with $\mathcal{V}_0 =span \{Y_1,Y_2\}$ and $\mathcal{V}_1 = span \{Y_3\}$ where we have the following nontrivial bracket operations, $[Y_1,Y_2]=Y_3, [Y_2,Y_3]=Y_2$.  Equation \eqref{rtrep} in section \ref{model} give a particular presentation of $\mathcal{RT}$.  We note that $\mathcal{RT}$ is not a Carnot group as it is not nilpotent.  
\end{exm}
Given a submanifold $S \subset \mathcal{G}$, the sub-Riemannian
geometry of $S$ is determined by the {\em horizontal normal} to $S$, $N_0$,
which is simply the projection of the Riemannian normal, $N$, to
$\mathcal{V}_0$, the first layer of the grading.  Explicitly, if
$\{X_1,X_2,X_3\}$ is a left invariant orthonormal basis of
$T\mathcal{G}$ with $\mathcal{V}_0 = span \{X_1,X_2\}$ and $S$ is
given as a level set $\phi=0$ then,
\begin{equation*}
\begin{split} N_0 = \text{proj}_{\mathcal{V}_0} N\\
& = proj_{\mathcal{V}_0} ((X_1 \phi) \;
X_1 + (X_2 \; \phi)  X_2 + (X_3 \phi) \; X_3) =  (X_1 \phi) \;
X_1 + (X_2 \phi) \; X_2
\end{split}
\end{equation*}

We also define the {\em unit horizontal normal}:
\[\nu = \frac{N_0}{<N_0,N_0>^\frac{1}{2}}\]

Minimal surfaces in Carnot-Carath\'eodory spaces have been
investigated in a number of settings\cite{CH,CHMY,DGN,GarNh,GP,Pauls:minimal,Pauls:Obstr,Pauls:cmcgroup}.  In particular, Danielli, Garofalo and Nheieu\cite{DGN} show that a $C^2$ hypersurface in a Carnot
group $G$ satisfy the following minimal surface equation:
\[ div_0 \nu =0\]
where $div_0$ is the horizontal divergence operator on $\mathcal{V}$.
As shown in the next two sections, this equation also characterizes
minimal surfaces in the class of groups described above which satisfy
assumption 1.

 \section{Minimal surfaces in $\mathcal{RT}$}\setS{PH}
 

\bgL{Pseudo}
If $\mathcal{G}$ satisfies assumption 1, i.e. \bgE{BC} [ \mathcal{V}_0 , \mathcal{V}_1] \subset \mathcal{V}_0\enE
then there exists  a global strictly pseudoconvex pseudohermitian structure $(\eta,J,\mathcal{G})$ with the following properties
\begin{itemize}
\item $H=\text{span }\mathcal{V}_0$ is the contact distribution for $\eta$.
\item $V= \text{span }\mathcal{V}_1$ is spanned by the characteristic vector field $T$ for $\eta$.
\item The Levi metric agrees with $g$ on $H$  and is conformal with constant scaling factor to $g$ on $V$. 
\end{itemize}
\enL

\pf Let $X_1$ and $X_2$ be a left invariant orthonormal frame for
$\mathcal{V}_0$ and set $T= [X_1,X_2] \in \mathcal{V}_1$. By left
invariance and the bracket generating property of $\mathcal{V}_0$, we
see that $X_1,X_2,T$ form a global orthonormal frame for
$T\mathcal{G}$. Set $\eta$ to be the dual $1$-form for $T$ with
respect to this frame. Then clearly $H= \text{Ker } \eta$ and strict
pseudoconvexity is immediate as $H$ bracket generates at $1$-step. The
additional bracket condition \rfE{BC} implies that $T \lrcorner
d\eta=0$. Thus the first two properties hold automatically regardless
of which complex structure $J$ is chosen for $H$.  

Define $J:\mathcal{V}_0 \to \mathcal{V}_0$ by $JX_1 = -X_2$  and $JX_2 =X_1$. Next we extend $J$ to all of $T\mathcal{G}$ by setting $JT=0$ and declaring $J$ to be linear over \rn{}. The Levi metric defined by
\[ h(X,Y) = d\eta(X, JY) + \eta(X)\eta(Y)\]
then clearly is compatible with $g$ in the required fashion.

\epf

In this setting, we can now employ the techniques of pseudohermitian
and CR geometry.  Our key tool is the existence of a canonical
connection $\nabla$, derived independently by Webster\cite{Webster}
and Tanaka\cite{Tanaka},  for any strictly pseudoconvex
pseudohermitian geometry. The defining properties of the connection
are as follows: 
\begin{itemize}
\item $H$, $T$, $\eta$, $d\eta$ and $J$ are all parallel.
\item $\text{Tor}(X,Y)=d\eta(X,Y)T$ for $X,Y \in H$.
\item $\text{Tor}(JX,T)=J\text{Tor}(X,T) \in H$ for $X \in H$.
\end{itemize}
Computations using this connection are most easily conducted in the
moving frame approach of Cartan, adapted to this setting by Webster
\cite{Webster}. For this technique, we first complexify the contact
distribution $H$ and define the space of $(1,0)$ and $(0,1)$ vector
fields to be the $+i$ and $-i$ eigenspace of $J$ respectively. The
$(1,0)$ vector fields are then spanned by 
\[ Z = X_2 -iX_1.\]
The vector fields $Z, \overline{Z}$ and $T$ then form an orthonormal
(complex) frame for the complexified  tangent space. The dual frame
will  be denoted by $\zeta, \overline{\zeta}$  and $\eta$. We
introduce the connection form $\omega$ via the identity 
\[ \nabla Z  = \omega \otimes Z.\]
With respect to our frame the Webster-Tanaka connection and Levi metric can be uniquely computed from the following equations\cite{Webster}:
\begin{itemize}
\item $d\eta = ih \zeta \wedge \overline{\zeta}$.
\item $dh = \omega h + h \overline{\omega}$.
\item $d\zeta = \zeta \wedge \omega + \eta \wedge \tau$.
\item $\tau =0$ mod $\overline{\zeta}$.
\end{itemize}
The $1$-form $\tau$ is known as the torsion form. 

This connection proves well adapted to many geometric problems. For the study of horizontally minimal surfaces we have the following theorem.
\bgT{Minimal} 
Suppose $S$ is a non-characteristic surface patch  in a Carnot group satisfying our structure conditions. The horizontal minimal surface equation for $S$ 
\bgE{HME} \text{div } \nu =0 \enE
can be written as
\[ \nabla_{J\nu} J\nu  =0 \]
where $\nu$ is the horizontal unit normal. Thus if $S$ satisfies $\rfE{HME}$ then $S$ is ruled by horizontal $\nabla$-geodesics.
\enT

\pf The volume form for the Levi metric is given by $dV = \eta \wedge
d\eta$. By the defining properties of the pseudohermitian structure
this is a constant multiple of the Riemannian volume form.  From this
we immediately see that $dV$ is parallel for $\nabla$ and the
divergence operator for $dV$ agrees with the Riemannian
divergencve. Further it f ollows that $ T \lrcorner \eta \wedge d\eta = d\eta$ and so
\[ \text{div }T =0.\]
Now $\nu$, $J\nu$ and $T$ form a local orthonormal frame for $T\mathcal{G}$. A standard formula in Riemannian geometry (see for example Kobayashi\cite{Kobayashi}) then yields
\[ \text{div } X =  \text{trace}( \nabla X + \text{Tor}(X,\cdot) ).\] 
If $X$ is horizontal then the second part of the trace formula vanishes identically by the defining properties of the Webster-Tanaka connection. Using our particular choice of frame we then see that
\begin{align*} 
\text{div } \nu &= \aip{\nabla_{\nu} \nu}{\nu}{} + \aip{\nabla_{J\nu} \nu}{J\nu}{} + \aip{\nabla_T \nu}{\nu}{}\\
&= -\aip{\nu}{\nabla_{J\nu} J\nu}{}
\end{align*}
as the second and third terms vanish because $H$ is parallel and the connection is metric respectively. Thus on a non-characteristic, horizontally minimal surface patch we have 
\bgE{JN}
\nabla_{J\nu} J\nu =0
\enE
everywhere. The integral curves of $J\nu$ are therefore $\nabla$-geodesics. But $J\nu$ spans the intersection of $TS$ with $H$. Thus the integral curves of $J\nu$ foliate $S$.
 
\epf

\bgR{CHMYnote} The divergence form of the mimimal surface equation in
$\mathcal{RT}$ was first shown by Citti and Sarti\cite{CS} (see section 2.9  proposition 3.1 of that paper) but is also
a consequence of the more general psuedoherimitian framework of Cheng,
Huang, Malchiodi and Yang\cite{CHMY}.  We also note that a version of
this theorem, showing that smooth minimal surfaces are ruled, was
first shown in section 2 of
Cheng, Huang, Malchiodi and Yang,\cite{CHMY} again in the more general
context of psuedoherimetian manifolds.  We include the proof here for completeness and because
it facilitates the computations below.
\enR

We shall now apply these techniques to the special case to the roto-translation group $\mathcal{RT}$. Here the underlying manifold is $\rn{2} \times \sn{1}$ and $\mathcal{V}_0$ is defined by setting
\[ X_1 = \cos \theta \pd{}{x} + \sin \theta \pd{}{y}, \quad X_2 = \pd{}{\theta}\]
and declaring them to be a left-invariant, orthonormal frame for a distribution $H$. The  Riemannian structure by defining the transverse vector field,
 \[T = [X_1,X_2] = \sin \theta \pd{}{x} - \cos \theta \pd{}{y}\]
 and declaring it to be unit length and orthogonal to $X_1$ and $X_2$. The remaining commutation relations can then be explicitly computed as
\[ [X_1,T]= 0, \quad [X_2,T] = X_1.\]
When we run through the construction of \rfL{Pseudo} we note that $g$ is exactly the Levi metric in this case. The contact form can be explicitly computed as \[ \eta = \sin \theta dx - \cos \theta dy\]
and the dual to the complex vector field $Z=X_2 -iX_1$ is
\[ \zeta = \frac{1}{2} \left( d\theta + i \cos \theta dx + i \sin \theta dy \right).\]
Straightforward computations then yield
\begin{align*}
 d\eta &= \cos \theta d\theta \wedge dx + \sin \theta d\theta \wedge dy = 2i \zeta \wedge \bt{\zeta}\\
 d\zeta &= \frac{i}{2} d\theta \wedge \left( -\sin \theta dx + \cos \theta dy\right) = -\frac{i}{2} \zeta \wedge \eta + \frac{i}{2} \eta \wedge \bt{\zeta}.
 \end{align*}
The first identity also follows from the fact that the pseudohermitian structure was explicitly constructed to ensure that $X_1$ and $X_2$ were orthonormal. Since $h=2$ we can immediately deduce from the 2nd Webster identity that the connection form $\omega$ is pure imaginary. Thus we can deduce that
\[ \omega = -\frac{i}{2} \eta, \quad \tau = \frac{i}{2} \bt{\zeta}.\]
This implies that for the frame $X_1,X_2,T$ the only non-trivial covariant derivatives are in the $T$ direction. By examining the real and imaginary parts of the equation 
\[ \nabla_{T} Z = - \frac{i}{2} Z\] we see $\nabla_{T} X_1 = \frac{1}{2}X_2$ and $\nabla_{T}X_2 = - \frac{1}{2} X_1$.

The horizontal $\nabla$-geodesics can be computed explicitly. Consider a curve $\gamma = (x,y,\theta)$. Thus 
\[ \dot{\gamma}  = (\dot{x}, \dot{y}, \dot{\theta}) = (\dot{x} \cos \theta + \dot{y} \sin \theta)  X_1 + \dot{\theta} X_2 + (\dot{x} \sin \theta - \dot{y} \cos \theta) T.\] 
Thus if $\gamma$ is a purely horizontal curve, we must have 
\[ \dot{x} \sin \theta - \dot{y} \cos \theta=0.\]
Under this assumption, $D_t \dot{\gamma} = 0$ if and only if both $\dot{\theta}$ and $\dot{x} \cos \theta + \dot{y} \sin \theta$ are constant. We can then solve the equation
\[ \begin{pmatrix} \sin \theta & - \cos \theta\\ \cos \theta & \sin \theta \end{pmatrix} \begin{pmatrix} \dot{x} \\ \dot{y} \end{pmatrix} = \begin{pmatrix} 0 \\ R_0 \end{pmatrix}\]
to obtain $\dot{x} = R_0 \cos \theta$, $\dot{y} =R_0\sin \theta$.  

\begin{itemize}
\item Case 1: $\dot{\theta} \ne 0$. Set $R = R_0 / \dot{\theta}$, then 
\begin{align*}
x &= x_c + R \sin \theta\\
y &= y_c - R \cos \theta\\
\theta& = \theta_0 + \dot{\theta} t.
\end{align*} 
Here $x_c = x_0 -R \sin \theta_0$, $y_c= y_0 + R \cos \theta_0$. 
\item Case 2: $\dot{\theta} =0$. Set $R=R_0$, then
\begin{align*}
x&= x_0 +R (\cos \theta_0) t\\
y&=y_0+R(\sin \theta_0) t\\
\theta &= \theta_0.
\end{align*}
\end{itemize}

In the sequel we shall refer to the horizontal $\nabla$-geodesics (and connected subsets of them) as rules.

\section{Missing data and amodal completion}\label{amodal}

We next turn to the problem of filling in missing image data.  Image
data may be missing for a number of reasons:  one object occludes
another, the existence of a ``blind spot'' in the retina or some other
physiological failure.  In terms of digital image processing, data
corruption, noise or object occlusion can {\em de facto} create a domain of
missing data.  

Using the roto-translation model for the
hypercolumn structure in V1 described above, image data is lifted to $\mathcal{RT}$
and missing data is filled by solving the minimal surface problem for
the given boundary data.\cite{CS}  Mathematically, if $D \subset \R^2$ is an
open domain where image data is missing, and $c \in \mathcal{RT}$ is
the image of $\partial D$ under the lift $\theta(x,y)$ defined in
section \ref{model}, then we wish to find a minimal surface $\Sigma$
so that $\partial \Sigma = c$.  Moreover, as discussed in the previous two
sections, a $C^2$ surface meeting these requirements must satisfy the equation:
\[div_0 \; \nu =0\]
where $\nu$ is the unit horizontal normal to $\Sigma$.  Moreover, by
theorem \rfT[PH]{Minimal}, $\Sigma$ must be ruled by horizontal
$\nabla$-geodesics.  For the balance of the paper, we will consider the following problem:
\vspace{.1in}
\begin{center}
\fbox{\parbox{5in}{
\noindent
{\bf Occlusion problem:}  Given a smooth curve, $c$, in $\mathcal{RT}$ which is the lift of the boundary of an open domain in $\R^2$, can we find a smooth minimal surface spanning $c$ which is ruled by $\nabla$-geodesics?
}}\end{center}

\vspace{.1in}

As referenced in the previous sections, this type of problem has been studied before in a number of sub-Riemannian settings.  In addition to the observation that minimal surfaces in some settings are ruled surfaces, there are a number of results further describing the nature of solutions to the minimal surface problem with Dirichlet boundary data.  Among these results, it is important to note that, at
least in the Heisenberg\cite{Pauls:Obstr} and Martinet-type spaces\cite{Cole}, there are
obstructions to the existence of smooth minimal spanning surfaces,
even if the spanned curve $c$ has arbitrarily nice behavior.
Moreover the second author demonstrates\cite{Pauls:minimal} that
solutions to the Dirichlet problem for ruled minimal surfaces need not
be unique (however, Cheng, Huang, Malchiodi and Yang\cite{CHMY} prove a uniqueness result for
surface subject to certain constraints on the characteristic locus).  Thus, as the roto-translation group is locally very much like the Heisenberg group,
we should expect to see issues with both existence and uniqueness. In light of this suggestive evidence, we present a list of conditions, each stronger than the next, concerning a smooth minimal spanning surface $\Sigma$:
\renewcommand{\theenumi}{\Roman{enumi}}
\begin{enumerate}
\item \label{conde} $\Sigma$ exists
\item \label{cond0} Condition \ref{conde} and any rule connecting two points of $c$ projects to a curve in the interior of $D$
\item \label{condA}Condition \ref{cond0} and the projection of $\Sigma$ to $D$ is surjective.
\item \label{condB} Condition \ref{condA} and $\Sigma$ is a graph over $D$
\item \label{condC} Condition \ref{condB} and $\Sigma$ is unique
\end{enumerate}
\renewcommand{\theenumi}{\arabic{enumi}}
We note that if condition \ref{cond0} is violated, the rules, upon projection,
would present potentially conflicting data for points exterior to $D$
while if condition \ref{condB} is violated, there would exist points
interior to $D$ with conflicting projected image data.  Thus, condition \ref{condA} is sufficient to guarantee the
existence of a completion of the image data (not just a spanning
surface in $\mathcal{RT}$) although there may be conflicting data
while \ref{condB} would provide a completion with no conflicting
data. However, as, {\em a priori} there may be multiple lifts, only
condition \ref{condC} would yield a unique completion of the image
data.   

The evidence cited above and the experimental evidence of
multiple simultaneous completions of image data in $\mathcal{RT}$\cite{CS} suggests that we should not expect to be able to satisfy the more stringent
requirements.  As we develop the machinery to construct such surfaces,
however, we will keep each of these conditions in mind.

\section{Ruled surfaces in $\mathcal{RT}$}\setS{PP} 
In this section, we begin the investigation of the existence and properties of minimal ruled spanning surfaces of curves in $\mathcal{RT}$.  As we require the solution to the occlusion problem to be a surface ruled by $\nabla$-geodesics, we first look at the set of points which can be connected to a given point by $\nabla$-geodesics.   With this in mind, we make the following definition:

\bgD{accessible}
For a point $p\in \mathcal{RT}$ we define the accessible set $\mathscr{A}(p)$ to be the collection of points that can be connected to $p$ by a single, horizontal $\nabla$-geodesic.
\enD

\bgL{Accessible}
Given a point $p=(x_0,y_0,\theta_0)$, the set of accessible points  is given by the implicit equation
\[ \frac{y-y_0}{x-x_0} = \tan \left( \frac{\theta+\theta_0}{2} \right).\]
\enL 

\pf
When the connecting $\nabla$-geodesic is a straight line this is immediate. The other case follows easily from the trigonometric identity 
\bgE{Trig}
\tan \left( \frac{\theta+\theta_0}{2} \right) = \dfrac{ \cos \theta_0 - \cos \theta}{\sin \theta - \sin \theta_0}.
\enE 
The proof of this identity is an easy exercise with the tangent half-angle formulas.

\epf

This provides a description of the accessible set of $p$: 
\bgL{Mobius}
Every accessible set $\mathscr{A}(p)$ is the image of an embedding of the M$\"\text{o}$bius strip into $\mathcal{RT}$.
\enL

\pf
We shall give two arguments for this result. One purely geometric, the other more analytic. 

From \rfL{Accessible} we note that each $\theta$-slice of $\mathscr{A}(p)$ projects to a straight line in the $(x,y)$-plane of gradient $\tan( (\theta+\theta_0)/2)$. As $\theta$ increases this line rotates spanning out a helicoid. However since we must identify $\theta=2 \pi$ with $\theta=0$ and the factor of $1/2$ inside the $\tan$ means that there is an orientation switch at the join. Thus $\mathscr{A}(p)$ is a non-orientable line over $\sn{1}$ which therefore must be diffeomorphic to a M$\"\text{o}$bius strip.

A more analytic approach is to consider the $\nabla$-exponential map at $p$ restricted to the horizontal distribution. From our explicit description of the horizontal geodesics passing through $p$ we note that if $p=(x_0,y_0,\theta_0)$ then
\bgE{ExpC} \begin{split} \text{exp}_p(aX_1+bX_2) =& \big(x_0 + a/b ( \sin(\theta_0 +b) - \sin \theta_0 ), \\
& \qquad y_0 + a/b(\cos \theta_0 - \cos (\theta_0+b) ) , \theta_0 + b\big)\end{split} \enE
at least when $b \ne 0$. When $b=0$ we instead get
\bgE{ExpL}  \text{exp}_p(aX_1) = (x_0 + a \cos \theta_0, y_0+a\sin \theta_0, \theta_0).\enE
If $\text{exp}(a,b)=\text{exp}(a\upp,b\upp)$ we must therefore have that $b = b\upp +2k\pi$ (with neither being $0$) and $a\upp/b\upp = a/b$. This later can be summarized as the points $(a,b)$ and $(a\upp,b\upp)$ must lie on the same line through the origin with $b = b\upp +2k\pi$. Therefore the exponential map is bijective from $\rn{} \times [ -\pi,\pi]$ to $\mathscr{A}(p)$ provided that the sides of the strip are identified via $(a,-\pi)\sim (-a,\pi)$. This provides an explicit embedding of the M$\"\text{o}$bius strip into $\mathcal{RT}$ with image $\mathscr{A}(p)$.

\epf

For each point on a curve $\gamma$, $\mathscr{A}(\gamma(t))$ may contain
many points of $\gamma$ or very few.  Of most interest to the question of building
spanning surfaces are the points $\gamma(t)$ where
$\gamma \cap \mathscr{A}(\gamma(t))= \{\gamma(t)\}$ - i.e. the points that only
connect to themselves.  These points give constraints on the formation
of a minimal ruled spanning surface.  To help understand these points,
we make the following definition: 

\bgD{RestrictedAccess}
Given an embedded curve $\gamma$ and a point $p \in \gamma$ we define 
\[ \mathscr{A}(p,\gamma) = \gamma \cap \mathscr{A}(p),\] the points in $\gamma$ accessible to $p$. A point $p$ such that $\mathscr{A}(p,\gamma) = \{p\}$ is called an solitary point of $\gamma$. The solitary points of $\gamma$ will be denoted $\mathscr{I}(\gamma)$. A point $p=\gamma(t)$ such that $\dot{\gamma} \in H$ is called a Legendrian point of $\gamma$. The Legendrian points of $\gamma$ will be denoted $\mathscr{L}(\gamma)$. We also define the orthogonal points of $\gamma$, denoted $\mathscr{O}(\gamma)$ to be where $\dot{\gamma} \in \text{span}\{X_1,X_3\}$.
\enD

As seen in the definition, there are two types of solitary points, the
Legendrian points and the non-Legendrian points. We remark that for a
Legendrian point $p$, the candidate rule passing through $p$ is
tangent to the curve $\gamma$ and, as in the Heisenberg
group\cite{Pauls:Obstr}, one can use this as a starting place for
building a ruled minimal spanning surface. the non-orientability of
the accessible sets in the roto-translation case means that unlike for
the Heisenberg group we cannot deduce that all solitary points are
Legendrian. Indeed  non-Legendrian solitary points present more of a
problem as the candidate rules will be transverse at such a point.  To
investigate the structure of the set of solitary points further we
prove the following lemma.

\bgL{Solitary}
For any embedded curve $\gamma$, the set $\mathscr{I}(\gamma)-\mathscr{L}(\gamma)$ is open.
\enL

\pf
 If $p \in \mathscr{I}(\gamma)-\mathscr{L}(\gamma)$, then $\gamma$ intersects $\mathscr{A}(p)$ only at $p$ and does transversely. A small perturbation of the base point $p$ will cause a small perturbation of the M$\"\text{o}$bius strip embedding. As $\gamma$ is transverse to $\mathscr{A}(p)$, a small perturbation of $p$ cannot increase the number of nearby intersections. Away from $p$ the Euclidean distance of $\gamma$ from $\mathscr{A}(p)$ can be uniformly bounded below and so a small perturbation will not introduce any distant intersections. Thus curve points sufficiently near to $p$ will also be solitary. Clearly they will also be non-Legendrian. 
 
\epf

\bgL{SolnonL}  Suppose $\gamma$ is an embedded curve and $\gamma(0) \in \mathscr{I}(\gamma) - \mathscr{L}(\gamma)$.  Then the map
\begin{equation*}
\begin{split}
 \Theta:\sn{1}&\ra \sn{1}\\
t & \mapsto \theta(\gamma(t))
\end{split}
\end{equation*}
is surjective.
\enL
\pf This is a topological argument.  Suppose $\gamma$ is a curve in
$\mathcal{RT}$ that projects to the xy-plane bounding a domain $D$.
Further suppose there exists $\theta_0 \in \sn{1}$ not in the image of
$\Theta$.  Then, there exists a neighborhood of $\theta_0$, $N_0$, so
that $N_0 \cap Im(\Theta) = \emptyset$.  By \rfL{Mobius},
$\mathscr{A}(\gamma(0))$ is a M$\"\text{o}$bius strip.    However
$\mathscr{A}_0(\gamma(0)):=\mathscr{A}(\gamma(0)) \minus (\rn{2} \times N_0)$
is orientable. Since
$\gamma$, by assumption, cannot enter $\rn{2} \times N_0$ we see that $\gamma$ lies to one side of
$\mathscr{A}_0(\gamma(0))$.  However, since $\gamma(0) \in
\mathscr{I}(\gamma) - \mathscr{L}(\gamma)$, $\gamma$ must intersect
$\mathscr{A}_0(\gamma(0))$ transversely at $\gamma(0)$.  This is a
contradiction.  

\epf

For the occlusion problem we shall work exclusively with curves that occur as the boundary of a smooth, simply connected bounded region $D$ lifted by the contour direction field of an intensity function $I:\rn{2} \to \rn{}$. The boundary $\partial D$ is can be viewed as the image of an embedding $\beta:\sn{1} \to \rn{2}$. Away from critical points of $I$, we can define the lifting function $\theta: \sn{1} \to \sn{1}$ by
\[ \theta(t) = \arctan \left( - \dfrac{I_x \circ \beta(t) }{I_y \circ \beta(t)} \right)\]
where at each point we choose the branch of $\arctan$ which makes $\theta$ continuous.  We shall that an occlusion is \textit{non-degenerate} if the number of critical points of $I$ lying inside $\partial D$ is finite and $\theta$ can be extended continuously across each critical point. The occlusion is \textit{completely non-degenerate} if there are no critical points on the boundary. Using this function $\theta$ we construct the curve $\gamma = (\beta,\theta)$.  In addition we define the normal angle function for $\beta$ by
\[ \dfrac{ \beta\upp(t) }{|\beta\upp(t)|} =   (  - \sin \varphi_{\beta} (t), \cos \varphi_{\beta}  (t) ).\]
From these we construct the \textit{transversality function} for $\gamma$,
\[ Q (t) = \theta(t) - \varphi_{\beta}(t).\] 
From the definitions, it is clear that a non-critical point $\gamma(t)$ is Legendrian if and only if $Q(t) = \pi/2 + k \pi$, $k\in \zn{}$. Likewise $\gamma(t)$ is orthogonal if and only if $Q(t) = k \pi$.  We record this and another fact in a lemma:

\bgL{LegendrianExist}
For any lift  $\gamma$ associated to a completely non-degenerate occlusion problem, \begin{enumerate}
\item $\gamma(t)$ is Legendrian if and only if $Q(t)=\frac{\pi}{2} + k \pi$ for some integer $k$.
\item $\gamma(t)$ is orthogonal if and only if $Q(t)= k \pi$ for some integer $k$.
\item $\mathcal{L}(\gamma)$ is non-empty. 
\end{enumerate}
\enL

\pf
For the first item, we note that, by definition, $\gamma$ is Legendrian if $\gamma' \in H$.  Computing, we have
\[\langle \gamma'(t), X_3 \rangle =  \beta'(t) \cdot (\nabla I) = |\beta'(t)|(\cos(\theta(t)-\phi_\beta(t)))=|\beta'(t)|\cos(Q(t))\]
and so, if we assume $\gamma(t)$ is Legendrian, we have that, equivalently,  $ \beta'(t) \cdot \nabla I=0$ or $Q(t)=\theta(t)-\phi_\beta(t)=\frac{\pi}{2} + k \pi$.  The second item, concerning orthogonal points, follows in the same way. 

The last item follows immediately from the observation that  
\[ \int_{\partial D} \nabla I \cdot d \vec{r} = \int_0^{2\pi} |\nabla I(\beta(t)) | |\beta\upp(t)| \cos \ua = 0\]
where $\ua$ is the angle between $\nabla I$ and $\beta\upp$. Since the
first two terms of the integral are strictly positive we must have
$\cos \ua $ taking both positive and negative values. In particular,
this implies there are at least two points in $[0,2\pi)$ where $\cos(\alpha)
=0$.  By the computation at the outset of the proof, these two points
are Legendrian points. 

\epf

\bgR{angleform}
We note that in the proof, we provide a geometric interpretation of $Q$: it measures the angle between $\beta'$ and $\nabla I$.
\enR
\section{Occluded Disks}\setS{CL}

For computational reasons, it is useful to restrict attention to
curves $\gamma$ that are lifts of circles in $\rn{2}$ to the
rototranslation group. The lifts we are most interested in come from
the direction angles of the contours of an intensity plot which has an ambiguity associated with the choice of orientation.  Given a point $(x,y)$ and a
contour passing through this point at angle $\theta$, it is unclear
whether to lift it to $(x,y,\theta)$ or
$(x,y,\theta+\pi)$. Accordingly, for a point $p = (x,y,\theta) \in
\mathcal{RT}$ we shall define its conjugate point to be $\overline{p}=
(x,y,\theta+\pi)$ and frequently consider conjugate lifts $\gamma$ and
$\overline{\gamma}$ simultaneously.  
 
Any circular lift $\gamma$ can be expressed parametrically in standard form as 
\bgE{Standard} \gamma(t) = \big( x_0 + R \cos t, y_0+R \sin t, \theta(t) \big).\enE
When $\gamma$ is understood, we shall frequently refer to a point of $p \in \gamma$ simply by its parameter value with this parametrisation. With this parametrization understood we can simplify the transversality function $Q$ to 
\bgE{DefQ} Q(t) = \theta(t) - t.\enE

\bgL{CircleAccess}
For a circular lift  $\gamma=(x,y,\theta)$, the non-trivial part of $\mathscr{A}(\gamma(t),\gamma)$  is given implicitly by
\bgE{Eqn} Q(t)+Q(u)= (2k+1)\pi, \qquad  k \in \zn{}\enE
and $\mathscr{A}(\gamma(t),\overline{\gamma})$ is given implicitly by
\bgE{CEqn} Q(t)+Q(u)= 2k\pi, \qquad  k \in \zn{}.\enE
\enL

\pf
From \rfL[PP]{Accessible} we see that $\gamma(u) \in \mathscr{A}(t,\gamma)$ if and only if
\[ \dfrac{\sin u - \sin t}{\cos u -\cos t} = \tan \left(\frac{\theta(u)+\theta(t)}{2}\right).\]
Applying the trigonometric identity \rfE[PP]{Trig} we see that this is equivalent to
\[ \cot \left( \frac{u +t }{2} \right)= \tan \left(\frac{\theta(u)+\theta(t)}{2}\right).\]
The result then follows easily from standard arguments in trigonometry. A virtually identical arguments yields the second part also.

\epf

\bgL{Conjugate}
Conjugation twist-commutes with the exponential map in the sense that.
\[ \text{exp}_{ \overline{p}}(a,b)  = \overline { \text{exp}_p(-a,b) }.\]
\enL

\pf
This follows from direct computation from \rfE[PP]{ExpC} and \rfE[PP]{ExpL}.
\epf

\bgC{ConjugateAccess}
If $ \overline{q} \in \mathscr{A}(p)$ then $q \in \mathscr{A}(\overline{p})$. Furthermore the projections to \rn{2} of the connecting rules match precisely.
\enC

Therefore when connecting points obtained from lifting intensity plots we need only consider how points in $\gamma$ can be connected to either $\gamma$ or $\overline{\gamma}$. We shall write $t \sim u$ if either $\gamma(u)$ of $\overline{\gamma}(u)$ lies inside $\mathscr{A}(\gamma(t))$.

\bgL{Legendre} 
The Legendrian points of $\gamma$ occur precisely where
\[ Q(t) = \pi/2 + k \pi, \qquad k \in \zn{}.\] 
The orthogonal points of $\gamma$ occur precisely where
\[ Q(t) =  k \pi, \qquad k \in \zn{}.\] 
\enL

\pf Since $\gamma$ is parametrized by \rfE{Standard}, we can explicitly compute that 
\begin{align*}
\dot{\gamma} &=  -R \sin (t) \pd{}{x} + R \cos(t) \pd{}{y} + \dot{\theta} X_2\\
&= R \sin (\theta -t) X_1 +R \cos( \theta -t) X_3 + \dot{\theta} X_2.
\end{align*}
For $\gamma(t) \in \mathscr{L}(\gamma)$ it is then necessary and sufficient that $\cos(\theta-t)=0$. For $\gamma(t) \in \mathscr{O}(\gamma)$ the condition becomes $\sin(\theta - t) =0$. The result follows easily.

\epf

\bgC{Legendre}
$ \overline{\mathscr{L}(\gamma)} = \mathscr{L}(\overline{\gamma})$, $ \overline{\mathscr{O}(\gamma)} = \mathscr{O}(\overline{\gamma})$. 
\enC

\bgC{LocalIsolated}
On an implicit plot of all points $(t,u)$ such that $t \sim u$. Any transverse crossing of the leading diagonal $u=t$ occurs at either a Legendrian or an orthogonal point.
\enC

\bgC{ConnectedLegendre}
If $t \sim u$ and $\gamma(t) \in \mathscr{L}(\gamma)$ then $\gamma(u) \in \mathscr{L}(\gamma)$. Furthermore, the connecting rule projects to the same circle as $\gamma$.

If $t \sim u$ and $\gamma(t) \in \mathscr{O}(\gamma)$ then $\gamma(u) \in \mathscr{O}(\gamma)$. 
\enC

With circular lifts, if two points on $\gamma$ are known to be connectable then it is a straightforward matter to explicitly describe the connecting rule. If $(x_0,y_0,\theta_0) \sim (x_1,y_1,\theta_1)$ with $\theta_0 \ne \theta_1$ then the connecting rule must have the parametrisation $(x_c + R \sin \varphi, y_c - R\cos \varphi, \varphi)$. Thus we need only solve the matrix equation
\begin{equation}\label{mateqn}
\begin{pmatrix} 1 & 0 & \sin \theta_0\\ 0 & 1 & -\cos \theta_0 \\ 1 & 0 & \sin \theta_1\\ 0 & 1 & -\cos \theta_1 \end{pmatrix} \begin{pmatrix} x_c\\ y_c\\ R \end{pmatrix} = \begin{pmatrix} x_0 \\ y_0 \\ x_1\\ y_1 \end{pmatrix}
\end{equation}
 where the fact that the points are connectable guarantees the existence of a solution. Elementary methods yield that generically
\bgE{Rule} R= \frac{x_1 -x_0}{\sin \theta_1 - \sin \theta_0}, \quad x_c = x_0 - R \sin \theta_0, \quad y_c = y_0 +R \cos \theta_0.\enE
This of course yields two separate connecting rule segments depending on whether $\varphi$ transverses \sn{1} clockwise or anticlockwise.

\bgL{Rnot0} If $\gamma$ is the circular boundary of a completely nondegenerate occlusion and if $R(t)=R(\gamma(t))$, then $R \neq 0$.  
\enL
\pf We may assume, without loss of generality, that
$(x_0,y_0,\theta_0)=(0,0,0)$, $\gamma(0)=(x_0,y_0)$ and $\gamma$ is parameterized by arclength.  If $R$
is zero, then we must have that $(x_0,y_0,\theta_0)$ is connected to
itself.  Thus, there exists a finite speed parametrization $(c_1(t),c_2(t))$ so
that $(c_1(0),c_2(0))=(0,0)$ and $(\gamma_1(c_1(t)),\gamma_2(c_2(t)))$ is connected to
$\gamma(t)$ for $t$ close to $0$. Using \eqref{mateqn}, if we let
$R(t)$ be the radius of the circle connecting $(\gamma_1(c_1(t)),\gamma_2(c_2(t)))$ to
$\gamma(t)$, for a generic choice of $t$, we must have 
\[R(t) = \frac{\gamma_2(c_2(t)) -
  \gamma_2(t)}{\cos(\theta(c_2(t)))-\cos(\theta(t))}=-\frac{\gamma_1(c_1(t)) -
  \gamma_1(t)}{\sin(\theta(c_1(t)))-\sin(\theta(t))}\]
Assuming that $R \ra 0$ as $t \ra 0$,  then we must have that
\begin{equation*}
\begin{split}
R(t) &= \lim_{t \ra 0}
\frac{\dot{\gamma}_2(c_2(t))\dot{c}_2(t)-\dot{\gamma}_2(t)}{\sin(\theta(c_2(t)))\dot{\theta}(c_2(t))\dot{c}_2(t)-\sin(\theta(t))\dot{\theta}(t)}\\
& =\lim_{t \ra 0}\frac{\dot{\gamma}_1(c_1(t))\dot{c}_1(t)-\dot{\gamma}_1(t)}{\cos(\theta(c_1(t)))\dot{\theta}(c_1(t))\dot{c}_1(t)-\cos(\theta(t))\dot{\theta}(t)}=0
\end{split}
\end{equation*}
As $\sin,\cos,\dot{c}_i$ are bounded and $|\dot{\gamma}|=1$, we have that $\dot{\theta}\ra \pm \infty$ as $t \ra 0$.
However, direct computation shows that 
\[ \dot{\theta}(t) = -\frac{1}{|\nabla I|} \begin{pmatrix}
 - \frac{I_x}{|\nabla I|}\\\frac{I_y}{|\nabla I|} \end{pmatrix}
   \mathscr{H} \begin{pmatrix} \cos(t)\\ -\sin(t) \end{pmatrix} \]
where $\mathscr{H}$ is the Hessian of $I$.  If $I \in C^2$ and the
   occlusion is completely nondegenerate, then $\dot{\theta}$ must be
   bounded on $\gamma$.  Thus, $R$ cannot tend to zero. 
\epf

With these initial observations in place, we turn to the task of understanding when minimal spanning surfaces exist.  We will investigate such surfaces using the following blueprint (if possible):

\begin{enumerate}
\item Construct the function $Q$ as given in \rfE{DefQ}.  
\item Find all solitary Legendrian points using lemma \rfL{Legendre}.
\item Starting from a solitary Legendrian points, construct connections between points on $\gamma$ to other points on $\gamma$ using the implicit equation 
\[ Q(t)+Q(u)=(2k+1)\pi\]
and the matrix equation above.
\end{enumerate}
Subject to goal of satisfying conditions \ref{cond0} and \ref{condA}
from section \ref{amodal}.  

As we shall see, any one of these steps and/or goals may be violated.
One issue we address first is that it may not always be possible to
construct such a surface by connecting points of $\gamma$ to other
points of $\gamma$.  In this case, as discussed above, it is natural
to instead connect some (or all) points of $\gamma$ to points of
$\overline{\gamma}$.

\section{The index of $Q$} \setS{IX}
 
To construct a minimal spanning surface for $\gamma$ we choose the a
Legendrian point as a starting place and attempt to build a monotone
function $u(t)$ such that $u \sim t$ by following the branch of the
implicit plot until we reach the other Legendrian point. If we
implicitly differentiate either \rfE[CL]{Eqn} or \rfE[CL]{CEqn} with
respect to $t$, we see that 
\[ u\upp(t) = \dfrac{\theta\upp(t)-1}{1-\theta\upp(u)} = -
\dfrac{Q\upp(t)}{Q\upp(u)}.\]  
In particular, these implicit plots fail to be graphs over the
$t$-axis precisely when either $Q\upp(t)=0$ or $Q\upp(u)=0$. The
presence of such points is a necessary condition for obstructions to
the existence of a monotone function $u(t)$.  If we have intervals of
positive measure where $Q'(t)$ is strictly positive and others
where it is strictly negative, then we note that a spanning surface
will not be a graph over $D$ as the change in sign forces the spanning
surface to backtrack locally, violating our goal \ref{condB}.  Generically, even if such
points occur we can still follow a branch of the implicit plot and
constuct connecting rules, however there will be non-uniqueness
issues. This approach will only completely fail if the plot fails to
be an embedded curve, which occurs only at points $(t,u)$ where
$Q\upp(t)=0=Q\upp(u)$.  

If we can construct a monotone function $u(t)$, then the collection of connecting rules will form a minimal spanning
surface for the lift $\gamma$.  In practice, there are several types
of obstructions to this method. 

Of primary importance is the degree of the map $Q\colon\sn{1} \to \sn{1}$. Since connections are made by implicitly solving the equation
\bgE{IEqn} Q(t) +Q(u) = \pi \quad (\text{mod } 2\pi)\enE the size of
the set $\mathscr{A}(\gamma,t)$ is intimately related to $\text{deg
}Q$. For example if $\text{deg }Q\ne0$ then there cannot be any
non-Legendrian solitary points as non-zero degree implies $Q$ is
surjective. Likewise, it is clear that  
\[ \left| \mathcal{L}(\gamma) \right| \geq 2 \left| \text{deg }Q \right|.\] 
Of course, critical points of $Q$ add a seperate pathology which
influences the number of branches of the implicitly defined function $u$ of \rfE{IEqn}. 

For the occlusion problem, the degree of $Q$ is directly related to the critical point theory of the intensity function $I(x,y)$. If the occlusion is completely non-degenerate, then there is a well-defined map
\[ I_\theta = \dfrac{\nabla I}{|\nabla I|} \colon \partial D \cong \sn{1} \to \sn{1}.\] 
The continuous lift $\theta$ is then given by $I_\theta + \pi/2$. Hence 
\[ \text{deg }Q = \text{deg }I_\theta -1.\]

To explore this further, we shall suppose that $I(x,y)$ has at most one (nearby) critical point $p$ which if it exists is contained in the interior of the occluded region $D$. In this instance it follows from the definitions that \[ \text{deg } Q = \text{index}_p \nabla I -1\]
if $p$ exists. If there is no critical point in the interior then $I_\theta$ extends to a continuous function on the interior disc. Standard results in algebraic topology then imply that $\text{deg }I_\theta =0$ and so $\text{deg }Q=-1$. 

Next, we examine various possibilities for the occlusion problem.  We
will deal primarily with completely nondegenerate curves.

\subsection{Case 1:  $deg \; Q = -1$, $Q' \neq 0$}

\begin{figure}
\centering
\includegraphics[width=3.0in]{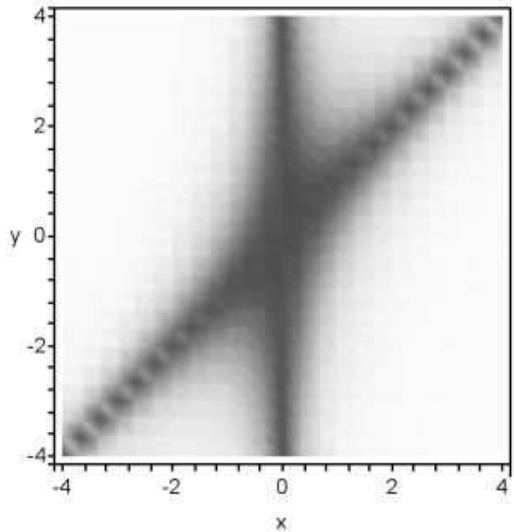}
\caption{ Intensity plot for $I=(1+x^2(x-y)^2)^{-1}$.}
\label{Fig:Cross}
\end{figure}

\begin{figure}
\centering
\mbox{\subfigure[$t\sim u$]{\includegraphics[width=2.0in]{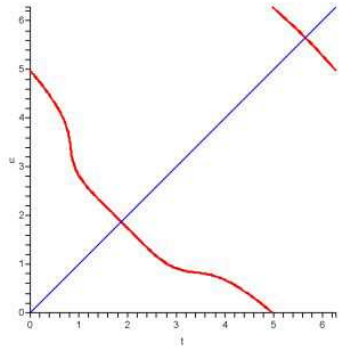}} \quad \subfigure[``$t$ vs $Q$'']{\includegraphics[width=2.0in]{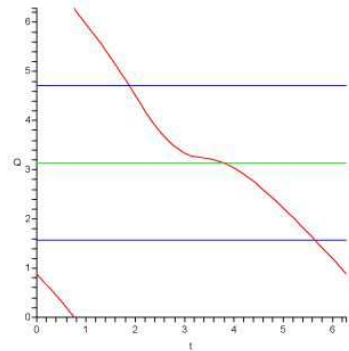}}}
\mbox{\subfigure[Contour Completion]{\includegraphics[width=2.0in]{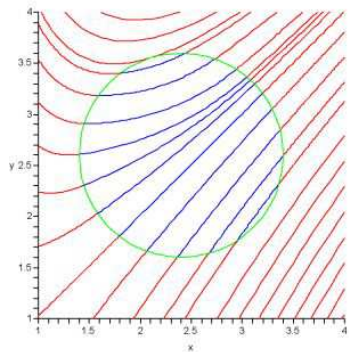}} \quad \subfigure[Minimal Spanning Surface]{\includegraphics[width=2.0in]{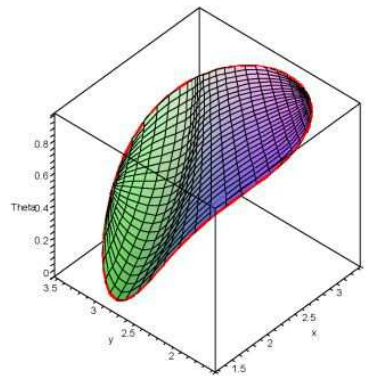}}}
\caption{ $I=(1+x^2(x-y)^2)^{-1}$, Centre=$(2.4,2.6)$, Radius=$1$, Connecting $\gamma$ to $\gamma$.}
\label{Fig:ISimple}
\end{figure}

As a basic example, we consider the following example: set the intensity function as $I(x,y)= (1+x^2(y-x)^2)^{-1}$ whose contour plot looks like an angled cross, see Figure \ref{Fig:Cross}. 
The region to be occluded is $(x-2.4)^2+(y-2.6)^2 <1$, which is over
just one branch of the cross.  While formally $I$ has critical points
everywhere along the lines $y=x$ and $x=0$, we can replace $I$ with
$x(y-x)$ without altering the underlying contour plot. With this
simplification, we have a completely non-degenerate occulsion with no
occluded critical points. Thus $\text{deg }Q=-1$ (and hence surjective
onto \sn{1}). The implicit plot of \rfE{IEqn} is shown in Figure
\ref{Fig:ISimple}. In this instance we see that $Q$ is one-to-one from
\sn{1} to \sn{1} and so there is only one branch of \rfE{IEqn}, which
spans the entire range.  Moreover, direct calculation shows that
$Q'(t) \neq 0$ for all $t\in [0,2\pi]$.  The boundary lift $\gamma$ has exactly $2$
Legendrian points, represented in Figure \ref{Fig:ISimple} (b) by the
intersection of the curve $Q(t)$ with the blue lines $Q=\pi/2$ and
$Q=3\pi/2$. The intersections with the lines $Q=0,2\pi$ and $Q=\pi$
correspond to the orthogonal points.  

Since Legendrian points can only connect to other Legendrian points
(by corollary \rfC[CL]{ConnectedLegendre}), we see from the implicit plot of \rfE{IEqn}
that the Legendrian points for $\gamma$ are solitary and correspond to
the intersections with the leading diagonal in Figure
\ref{Fig:ISimple} (a). Thus we can pick either and construct
connecting rules by tracing the sole branch of the implicit plot until
we reach the other Legendrian point. By symmetry every point on the
curve has now been connected to another and we can build a surface
ruled by the $\nabla$-geodesic segments that project into the interior
of the occluded region. The projections of these segments provide a
contour completion through the occluded region. In Figure
\ref{Fig:ISimple} (c) and (d)  we show the full contour completion
together with the associated minimal lift in the roto-translation
group. We note that it is possible to show that the surface constructed is a graph over the occluded region and hence satisfies condition \ref{condB}.  

\subsection{Case 2: $deg\; Q = -1$, $Q'$ has zeros}
\begin{figure}
\centering
\mbox{\subfigure[$t\sim u$]{\includegraphics[width=2.0in]{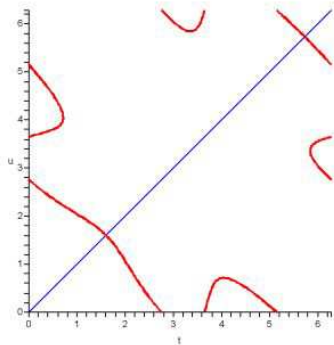}} \quad \subfigure[``$t$ vs $Q$'']{\includegraphics[width=2.0in]{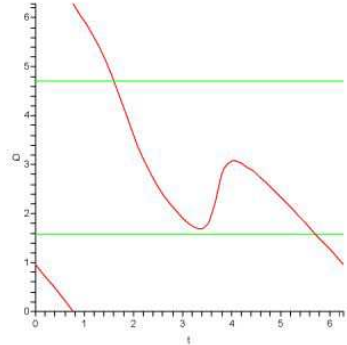}}}
\mbox{\subfigure[Contour Completion]{\includegraphics[width=2.0in]{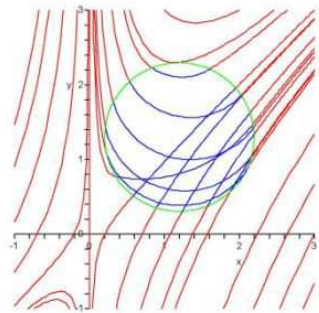}} \quad \subfigure[Minimal Spanning Surface]{\includegraphics[width=2.0in]{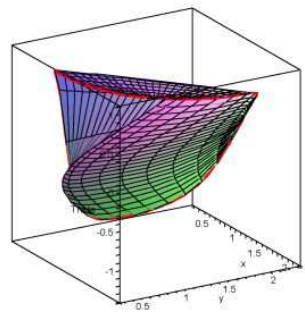}}}
\caption{ $I=(1+x^2(x-y)^2)^{-1}$, Centre=$(1.2,1,4)$, Radius=$1$, Connecting $\gamma$ to $\gamma$.}
\label{Fig:Back}
\end{figure}

The last example worked so well because not only was $\text{deg }Q=-1$,
but also because $Q$ was injective. This was because the directions of
the underlying contours were relatively uniform. If we move the
occluded region in closer to the centre of the cross, we lose this
uniformity and we find that $Q$ depsite having degree $-1$ is no
longer injective. Examining Figure \ref{Fig:Back} confirms that there
are points where $Q'=0$ (this can also be confirmed by direct calculation).  This is represented by
the failure  of the implicit plot of \rfE{IEqn} to be a graph over
either $u$ or $t$. When we follow the program laid out earlier for
constructing minimal spanning surfaces, we find that some points have
multiple connections. The surface then connects to some parts of the
curve as a ridge.  Note this surface satisfies condition \ref{cond0} but not \ref{condA}.  

\subsection{Case 3:  $deg \;Q=0$}
If $\text{deg Q}=0$ then there are possible obstructions to even local existence of spanning surfaces for $\gamma$. This phenomenon occurs due to the presence of non-Legendrian solitary points. Non-Legendrian solitary points do  not occur in isolation but as open sets by \rfL[PP]{Solitary}. Since $Q$ is continuous, the presence of non-Legendrian solitary points implies that the image of $Q$ is contained in a narrow (width $< \pi$) band. In particular, this implies that the condition $\text{deg }Q=0$ is necessary. In the situation of a single critcal point of $I$ being occluded, the index of $\nabla I$ must be $1$. In other words, we must be occluding a local maximum or minimum.

A simple example is to consider the circular lift
\bgE{PureSolitary} \gamma(t) = ( \cos t, \sin t, t ).\enE
Here $Q(t)=0$ everywhere and so every point is orthogonal with outward pointing orientation. The set of $(t,u)$ that satisfy \rfE{IEqn} is therefore empty. Every point is therefore solitary and non-Legendrian and there is no non-characteristic minimal surface that spans even a part of $\gamma$. 

\begin{figure}
\centering
\mbox{\subfigure[ $t \sim u$]{\includegraphics[width=1.8in]{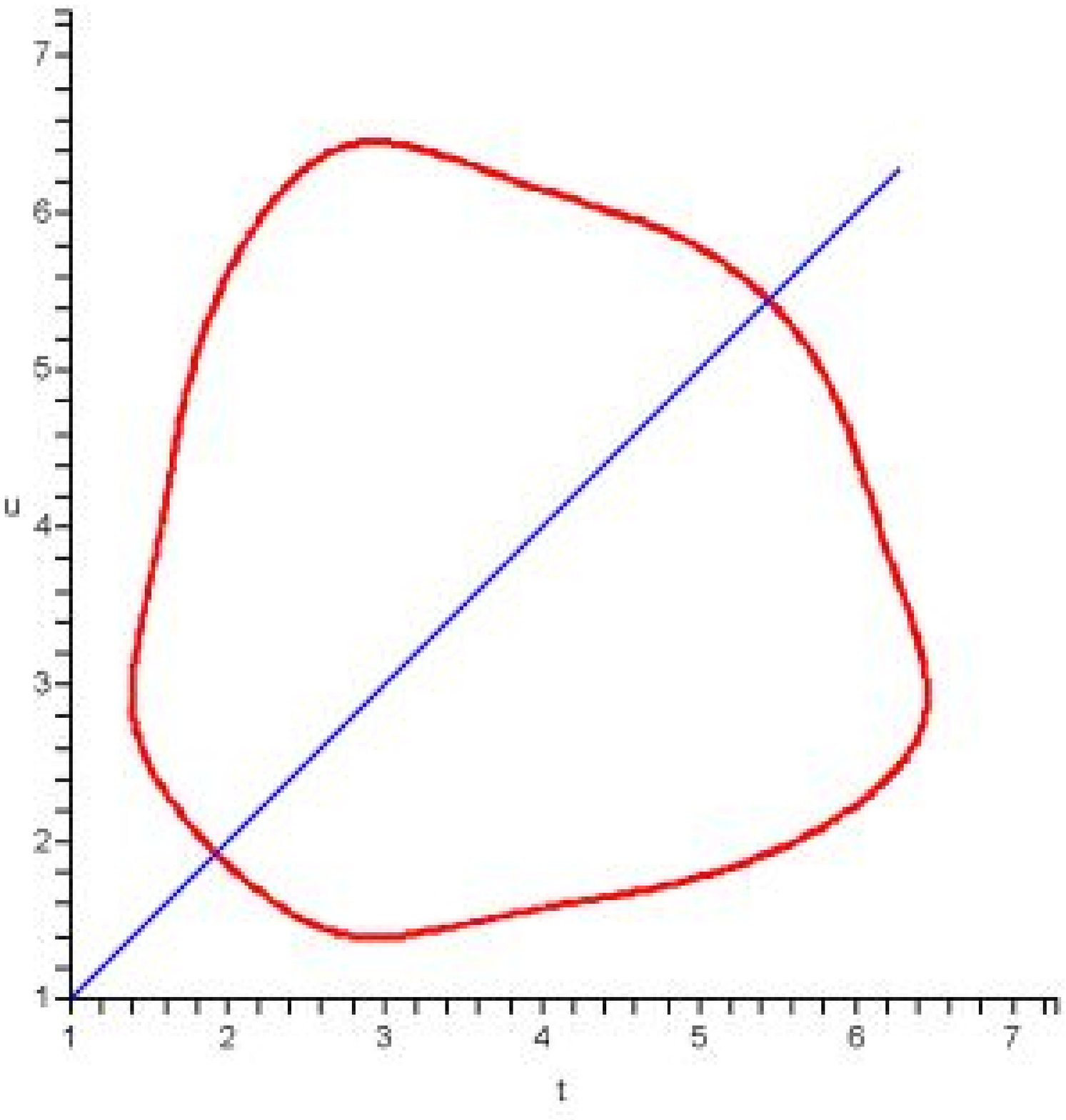}} \quad 
\subfigure[``$Q$ vs $t$'']{\includegraphics[width=1.8in]{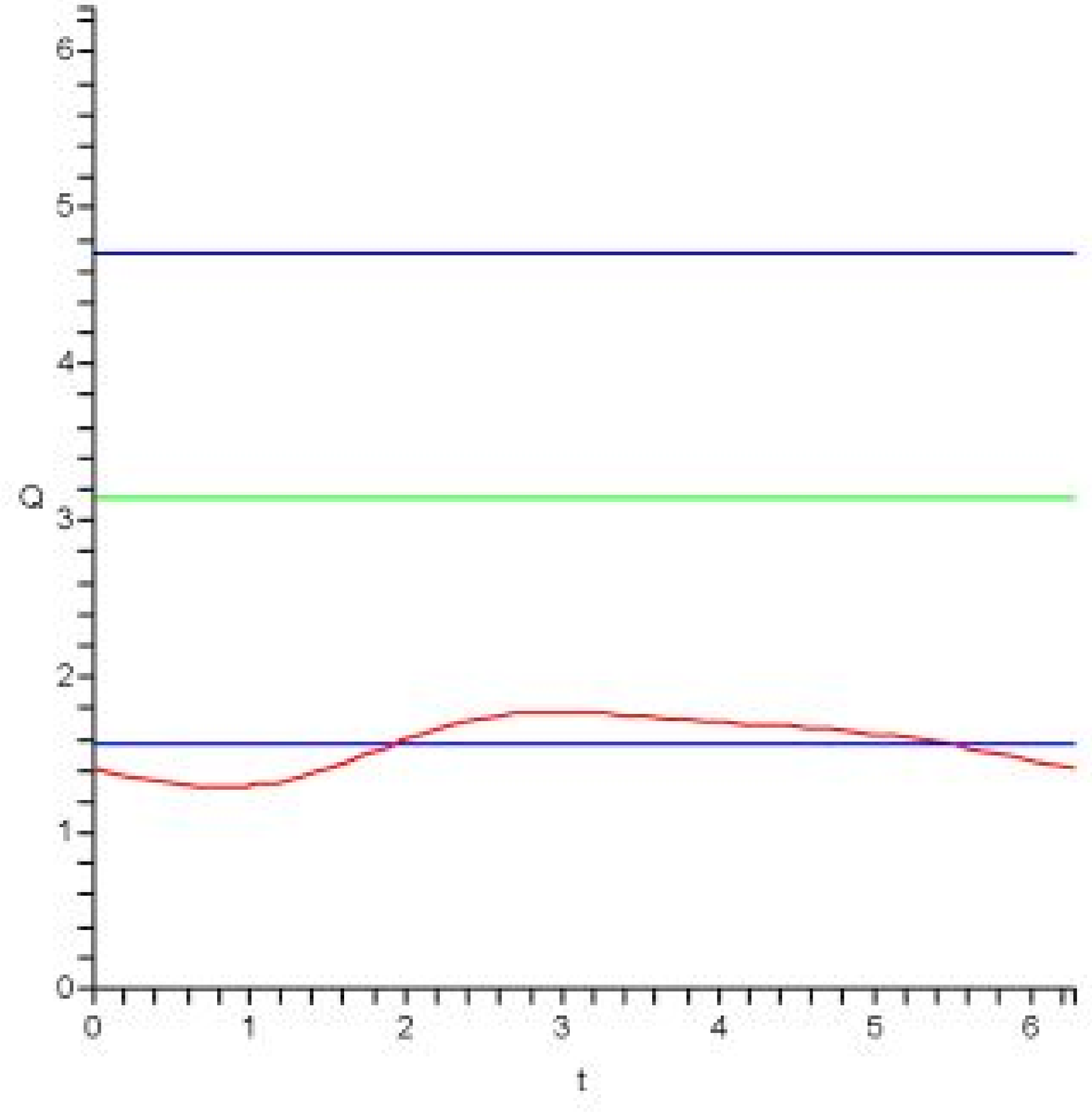}}}
\mbox{\subfigure[ Contour Completion]{\includegraphics[width=1.4in]{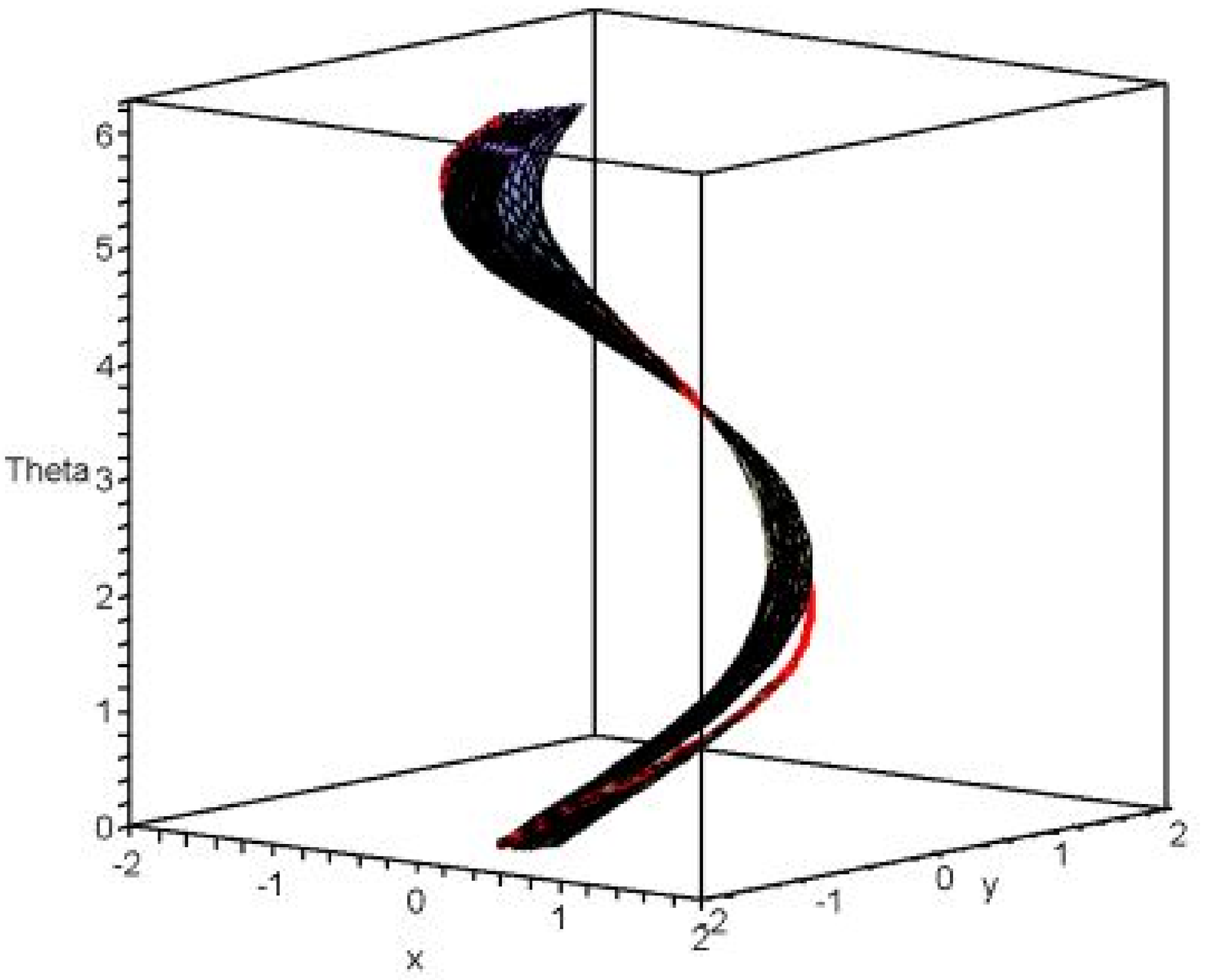}} \quad 
\subfigure[Close Up]{\includegraphics[width=1in]{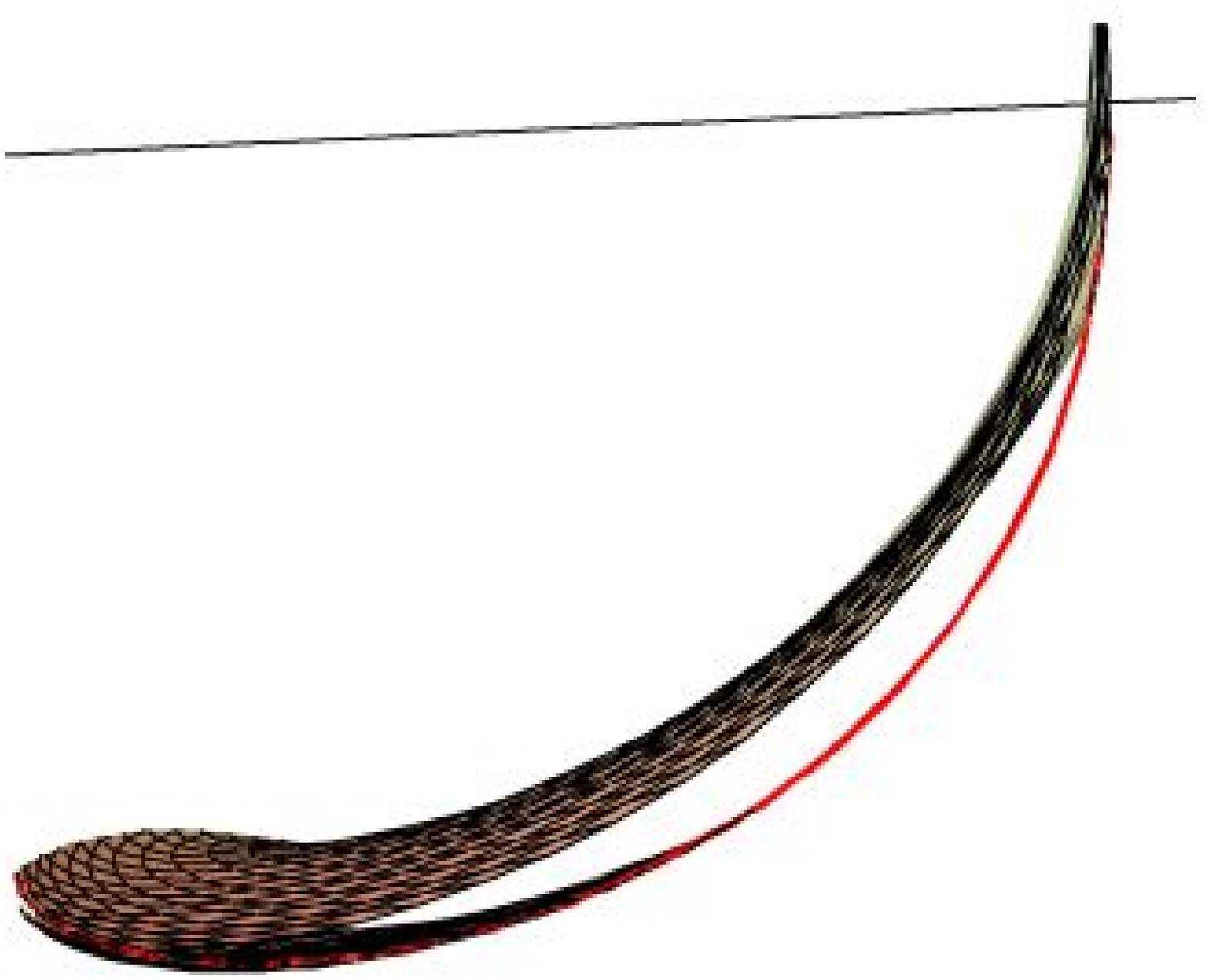}}
\quad 
\subfigure[Over View]{\includegraphics[width=1.4in]{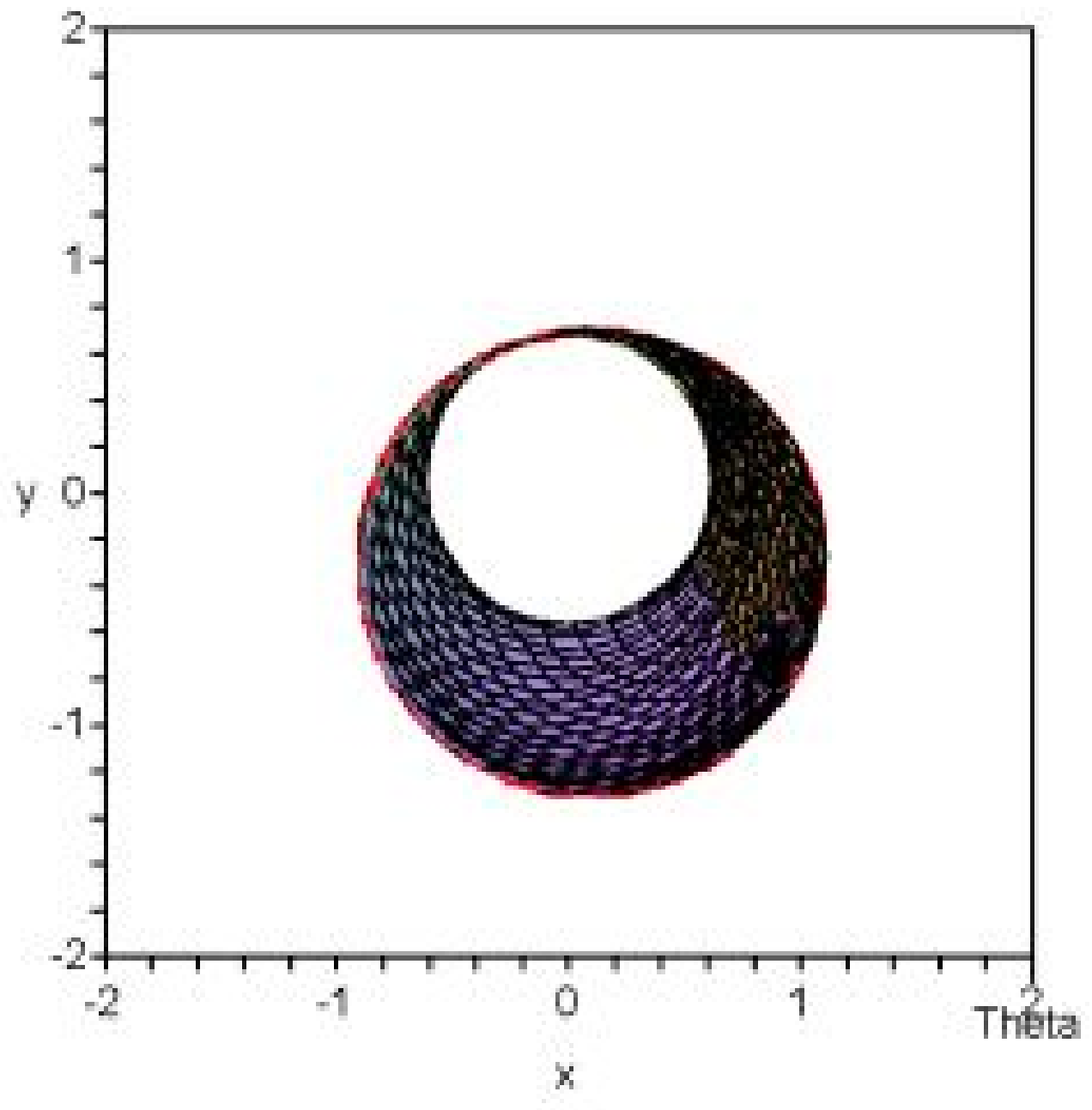}}}
\caption{ $I=(1+x^2+0.9y^2)^{-1}$, Centre=$(0.1,-0.3)$, Radius=$1$, Connecting $\gamma$ to ${\gamma}$.}
\label{Fig:Gap}
\end{figure}

It is clear that the lift \rfE{PureSolitary} cannot occur from an
occlusion problem as it would require the vector field $\nabla I$ to
be rotational and hence non-conservative. However gaps in the implicit
plot are characteristic of occluded maxima and minima, at least in the
absence of symmetry. See Figure \ref{Fig:Gap} for an explicit example
where $I=(1+x^2+0.9y^2)^{-1}$ and the circle occludes a local maximum
of the function.  In figure \ref{Fig:Gap} (a), we see that there are
two gaps where there are no connections between $\gamma(t)$ and any
other point on the curve.  The nature of the gap as part of the
minimal spanning surface is shown in the remaining graphs.  

\subsection{Case 4: $|deg \; Q| > 1$}
\begin{figure}
\centering
\mbox{\subfigure[$t\sim u$]{\includegraphics[width=2.0in]{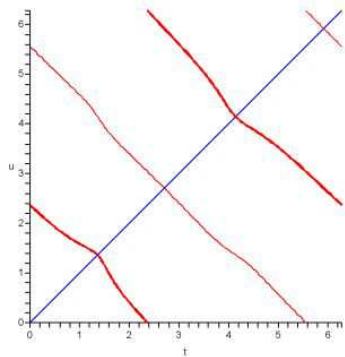}} \quad \subfigure[``$t$ vs $Q$'']{\includegraphics[width=2.0in]{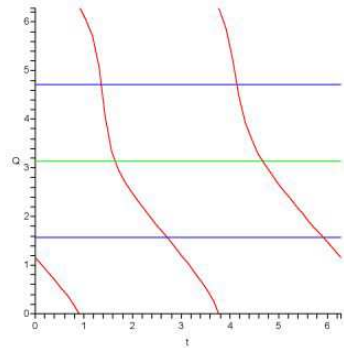}}}
\mbox{\subfigure[Contour Completion I]{\includegraphics[width=2.0in]{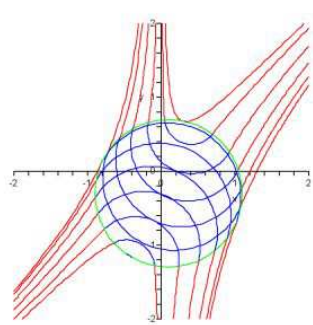}} \quad \subfigure[Contour Completion II]{\includegraphics[width=2.0in]{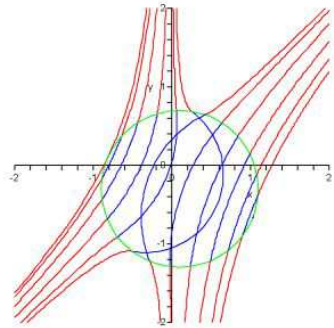}}}
\caption{ $I=(1+x^2(x-y)^2)^{-1}$, Centre=$(0.1,-0.3
)$, Radius=$1$, Connecting $\gamma$ to $\gamma$.}
\label{Fig:Double}
\end{figure}

When $Q$ has large degree, the phenomena of overlapping contours and
immersed, discontinuous spanning surfaces occurs naturally even when
we are considering only rules connecting $\gamma$ to itself. In Figure
\ref{Fig:Double} we return to the intensity function $I=
(1+x^2(x-y)^2)^{-1}$, but with the occluded region shifted to have
centre $(0.1,-0.3)$ and radius $1$. Since we are now occluding the
saddle point (of $x(y-x)$) at $(0,0)$, the degree of $Q$ is
$-2$. 

In this instance we see that $Q$ is now a monotone two-to-one function from \sn{1} to \sn{1}. This corresponds to there now being two branches of the implicit plot of \rfE{IEqn} in Figure \ref{Fig:Double} (a). For convenience of reference we shall refer to the highlighted branch as branch I and the other as branch II. In both contour completions there is overlap as Legendrian points are crossed while transversing the branches.

\subsection{Connecting $\gamma$ to $\overline{\gamma}$}
\begin{figure}
\centering
\mbox{\subfigure[$t \sim u$]{\includegraphics[width=1.8in]{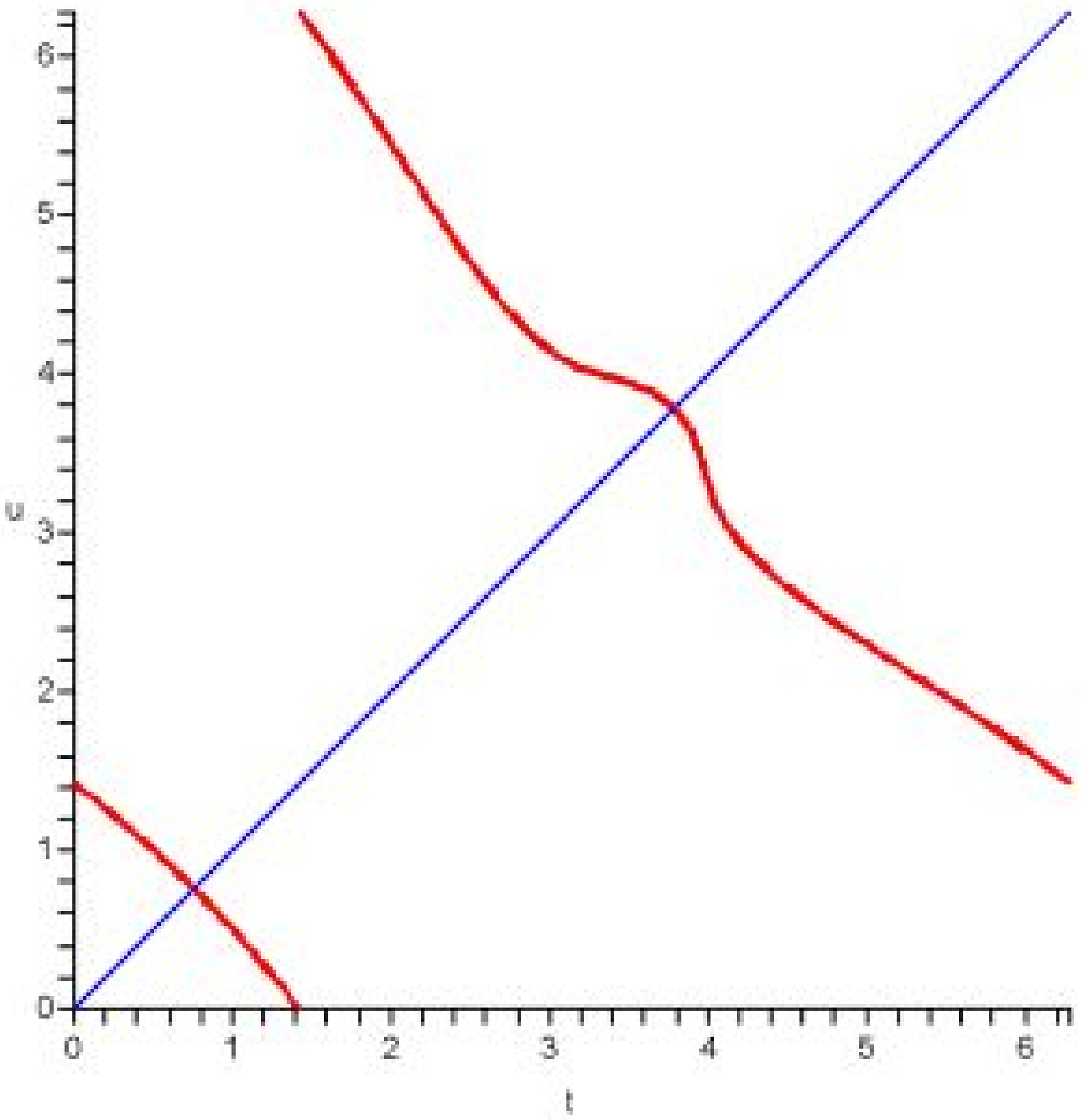}} \quad 
\subfigure[Contour Completion]{\includegraphics[width=1.8in]{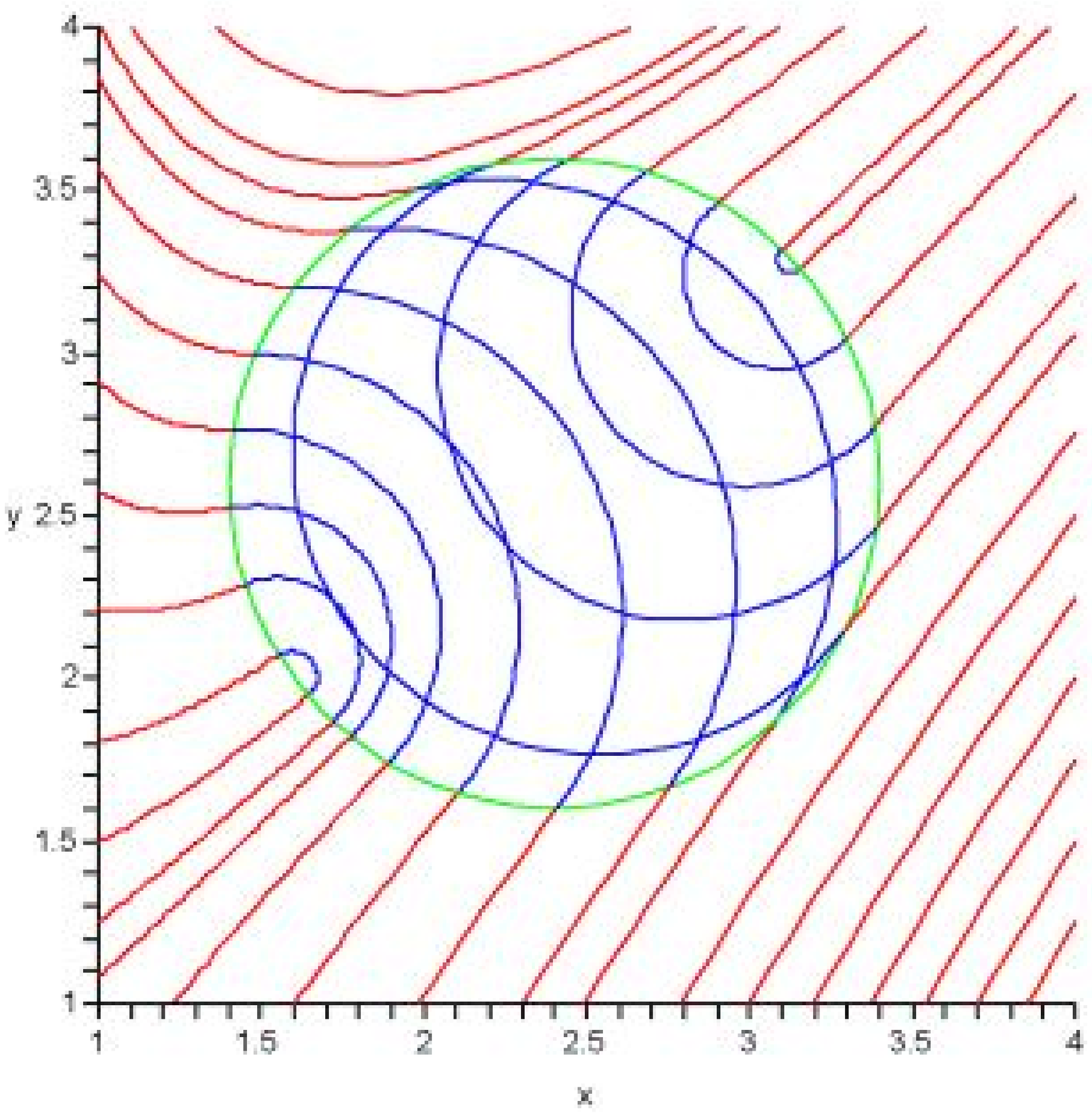}}}
\mbox{\subfigure[Minimal Spanning Surface I]{\includegraphics[width=2in]{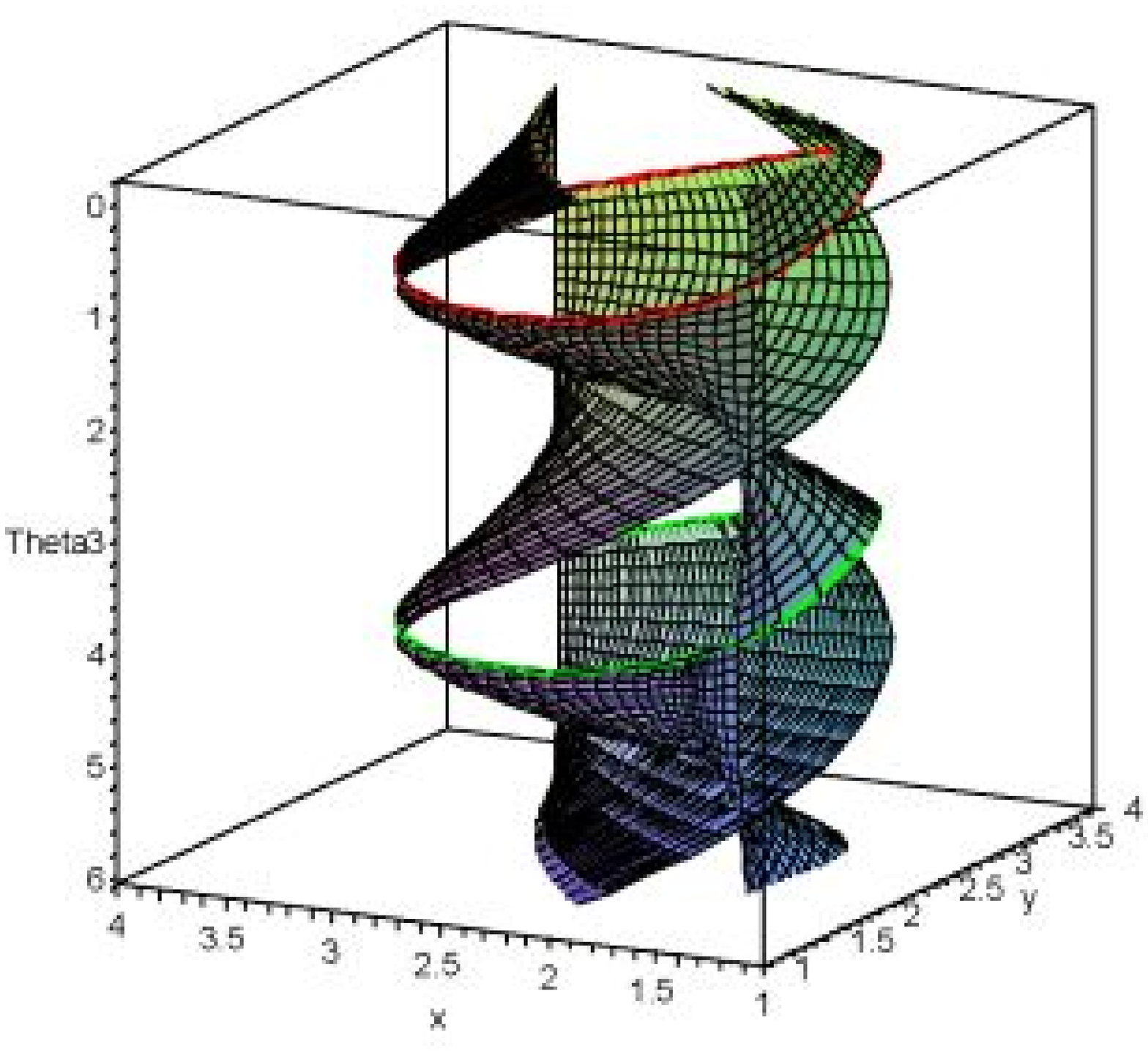}} \quad \subfigure[Minimal Spanning Surface II]{\includegraphics[width=2in]{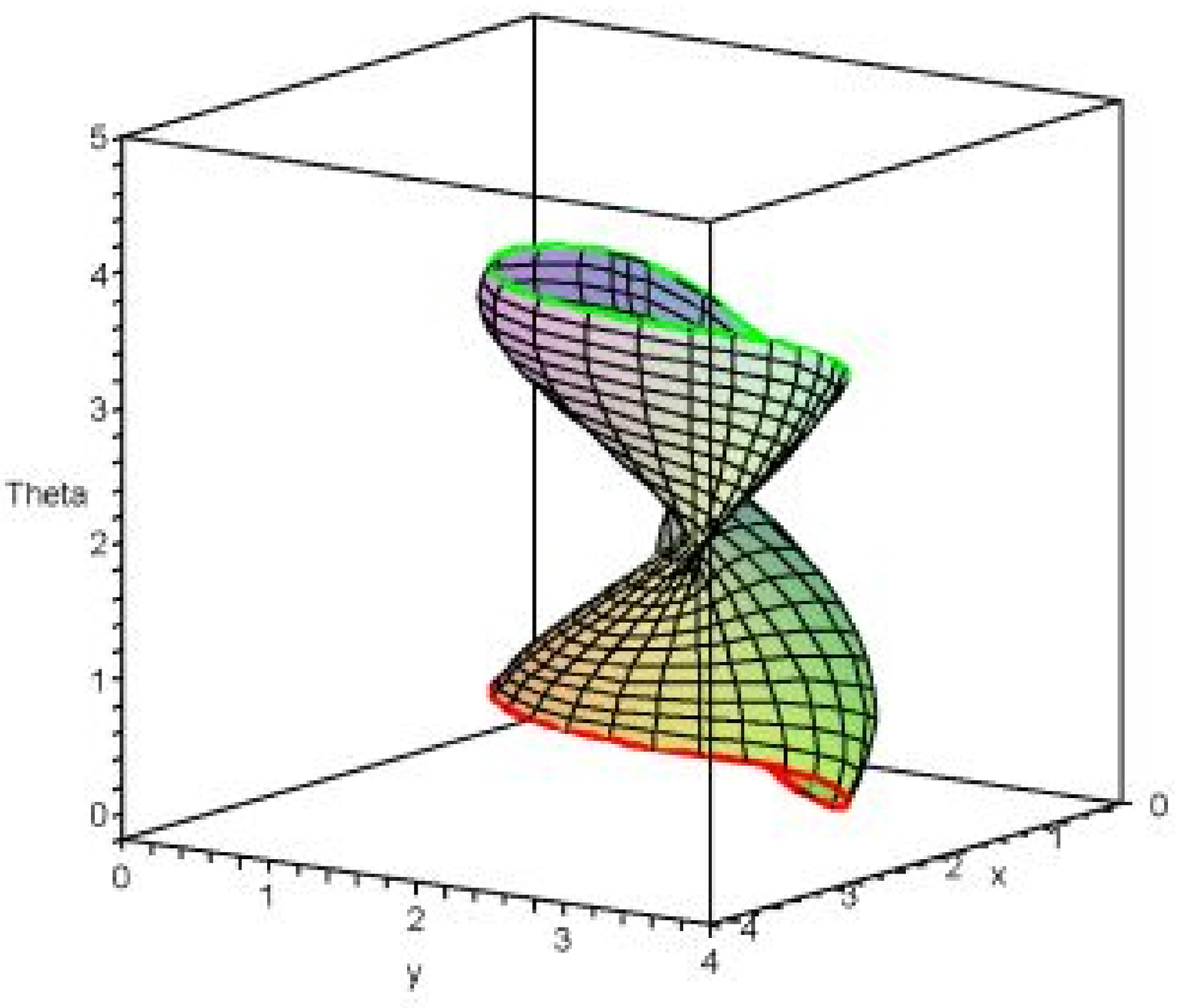}}}
\caption{ $I=(1+x^2(x-y)^2)^{-1}$, Centre=$(2.4,2.6)$, Radius=$1$, Connecting $\gamma$ to $\overline{\gamma}$.}
\label{Fig:SimpleC}
\end{figure}

For the problem of visual completion we must consider  also rules
connecting the lift $\gamma$ to its conjugate lift
$\overline{\gamma}$. To illustrate this we return to our original
example with the occluded region being the unit disc centred at
$(2.4,2.6)$. We follow the same basic program, but instead focus on
the orthogonal points as our start and finish
locations. Unfortunately, this introduces a pathology into the
construction of our contour completion and minimal spanning
surfaces. To progress from one orthogonal point to the other along a
branch of the implicit plot of \rfE[CL]{CEqn} it is necessary to pass
through a Legendrian point. The effect of this is to switch which
segment of the connecting rules projects to the interior of the
occluded region. As is seen in Figure \ref{Fig:SimpleC} (b) this
causes overlaps of the contour completion. If we look at the minimal
surface formed from these internal rule segments, in order to
completely span $\gamma$ we must continue further along the branch
until we return to the original orthogonal point.  This produces the
discontinuous self-intersecting surface of Figure \ref{Fig:SimpleC}
(c) . From the perspective of minimal spanning surfaces, it is more
natural to allow rule segments that project outside the occluded
region. As is shown in Figure \ref{Fig:SimpleC} (d), this yields a
smooth immersed surface between $\gamma$ and $\overline{\gamma}$, but
it still self-intersects. 

\section{A sufficient condition for effective disocclusion}\setS{thms}
The examples of the last section point to several general features of
the solutions to the occlusion problem.  We will now show that under some assumptions on the image function $I$, we
can guarantee the existence of a solution to the occlusion problem,
ie. a minimal spanning surface satisfying condition \ref{condA}.  To
do so, we need some preliminary lemmata.

\begin{figure}
\begin{center}
\mbox{\subfigure[$s(t)$ switchs
  signs]{\includegraphics[width=2in]{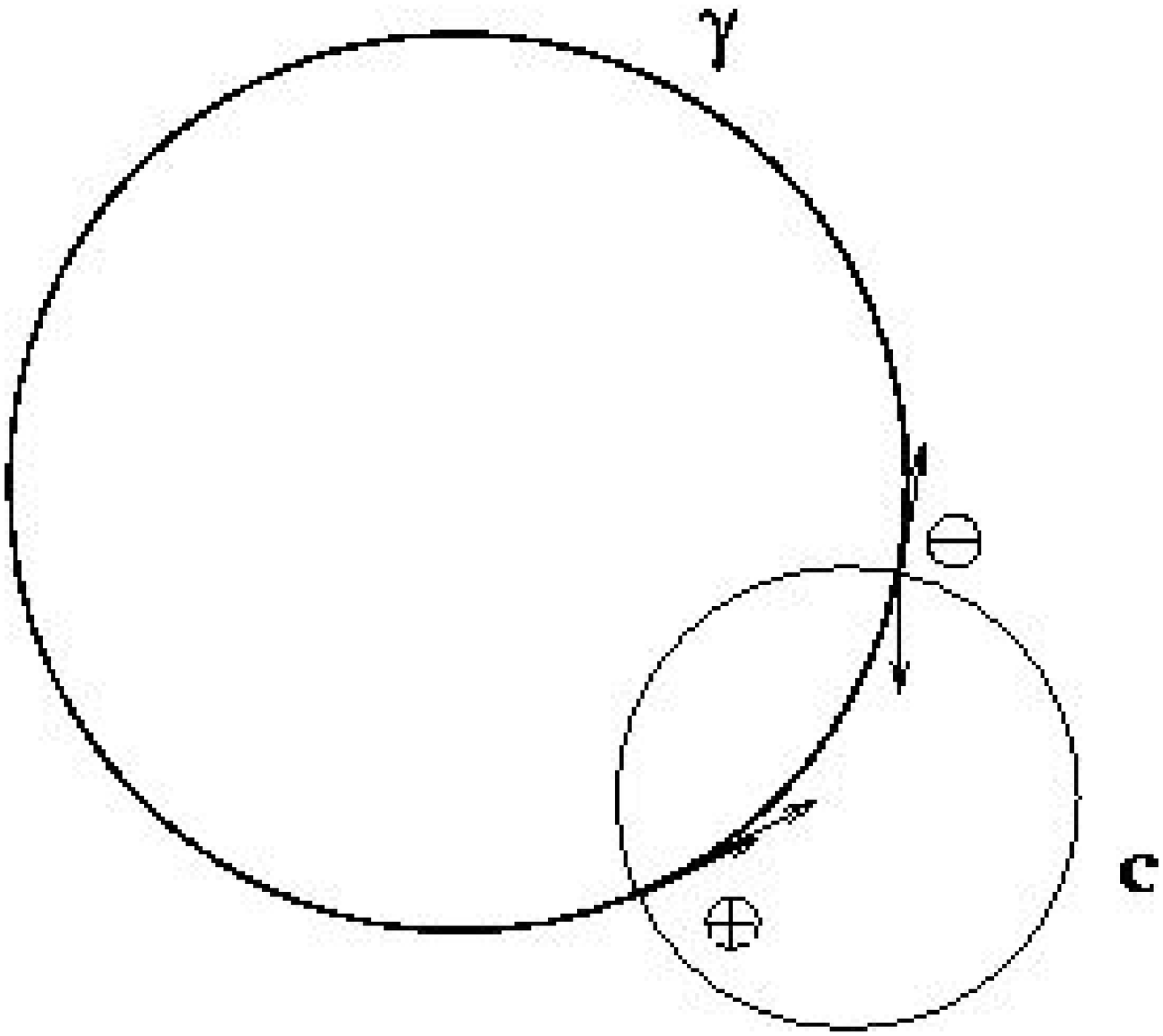} }\quad
  \subfigure[Two rules]{\includegraphics[width=2in]{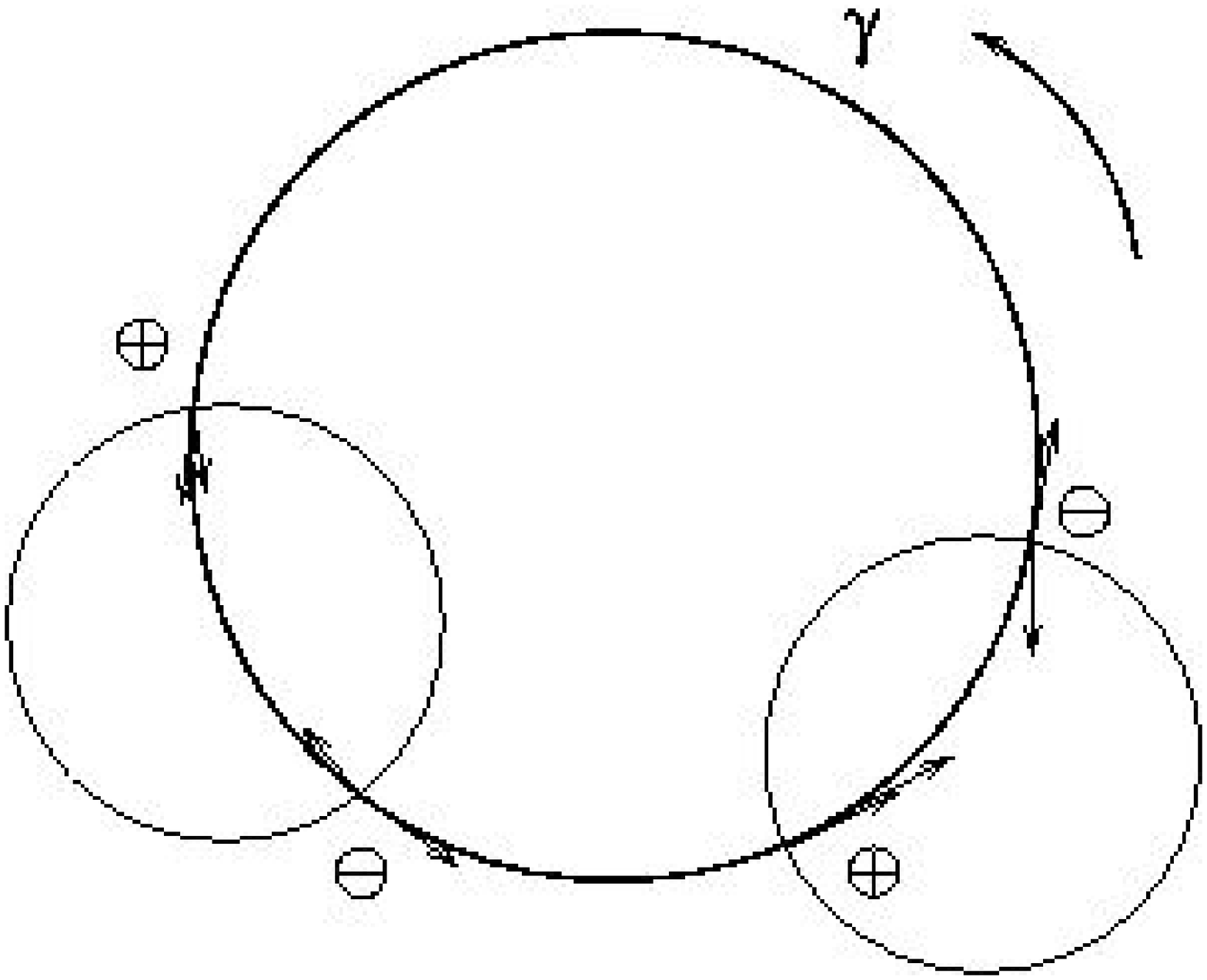}}}
\mbox{\subfigure[$|R|\ra\infty$]{\includegraphics[width=2in]{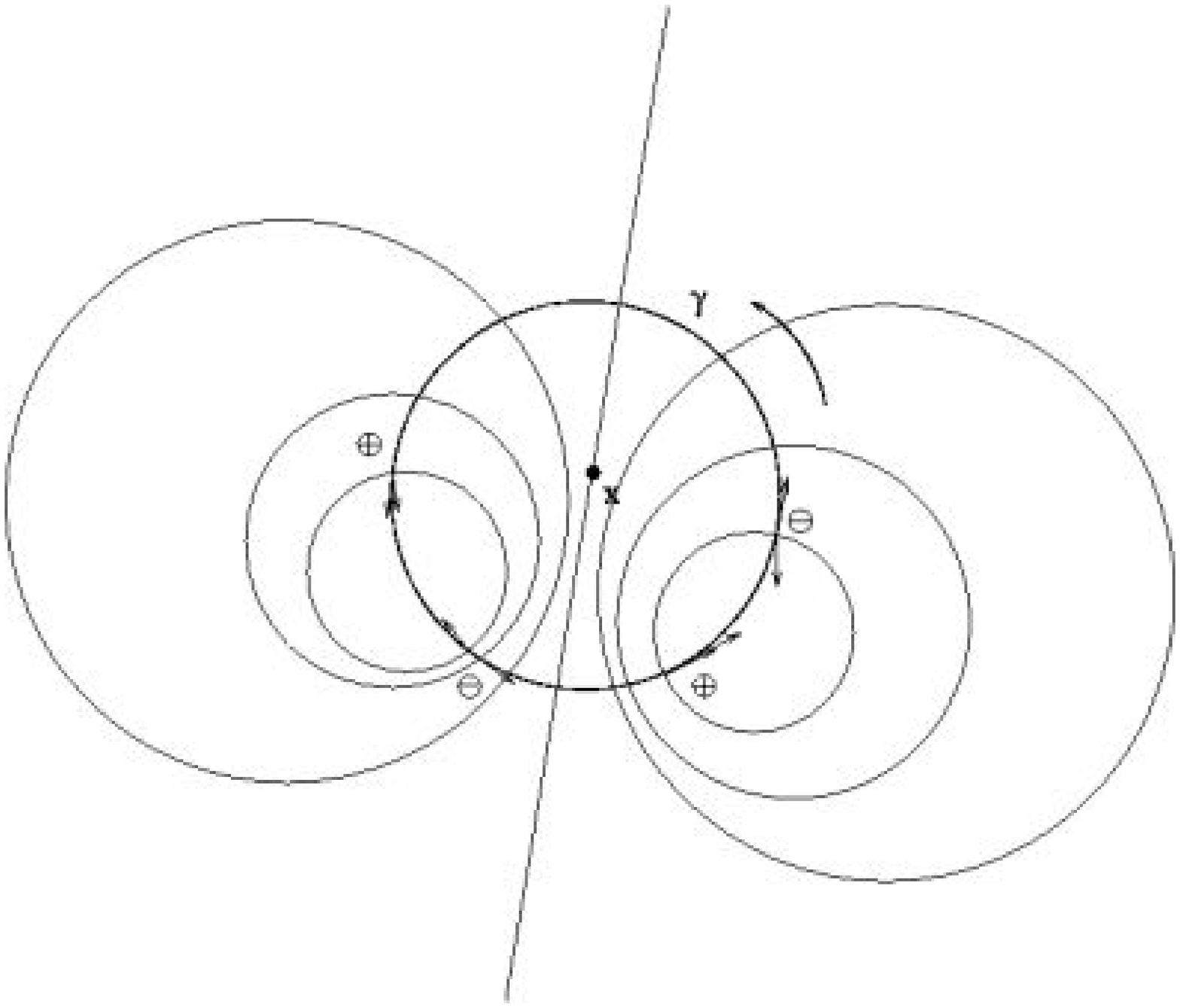}}}
\caption{Behavior of rules in a spanning surfaces}\label{Fig:switchsign}
\end{center}
\end{figure}

\bgL{surj1}  Let $I: \R^2 \ra \R$ be an intensity function of an
image with a completely nondegenerate occlusion given by a circular region $D$.  Suppose there
exists a minimal spanning surface $\Sigma$ of the occluded region associated to a monotone function $t \ra u(t)$ and so
that the projection of $\Sigma$ to $D$ is not surjective.  Then, one
of the rules of the minimal spanning surface is a circle that, when
projected to $\R^2$ lies entirely inside $D$ and is tangent to the
boundary of $D$ at a Legendrian point.
\enL 
\pf  We begin with some simple geometric observations.  First, if $\Sigma$ is composed of
circles of infinite radius (i.e. straight lines), the
projection is trivially surjective.  So, we may assume
there are some circles in $\Sigma$ that have finitie radii.  Second, let $\gamma(t)$ be the
circle that bounds the occluded region (as above), oriented counterclockwise, and let $c$ be a
circle of finite radius connecting two points on $\gamma$ that is the projection of a
rule of $\Sigma$.  If $\vec{n}$ is the inward
pointing normal to $c$, then let $s(t)=\cos(\alpha(t))$ where
$\alpha(t)$ is the angle (in $\R^2$) between $\dot{\gamma}$ and $\vec{n}$
at $\gamma(t)$.  Then, if $c$ connects $\gamma(t_1)$ to $\gamma(t_2)$,
$\text{sign}(s(t_1))=-\text{sign}(s(t_2))$ (see figure \ref{Fig:switchsign} (a)).  Moreover, in the
degenerate case where $t_1=t_2$, the circle $c$ is tangent to $\gamma$ (i.e. $\gamma(t_i)$ is Legendrian) and thus $s(t_i)=0$.  Since $u$ is monotone (and in particular, one to one) and $R,x_c,y_c$ are continuous in $t$, $s$ is well defined and, for $0 \le |R| < \infty$, $s$ is continuous.  At points where $|R|=\infty$, $s$ may have a jump discontinuity as $s$ will switch sign at such points.

By \rfL[PP]{LegendrianExist}, we have that there are at least two
Legendrian points and so, by the previous discussion, the function $s$ potenitally switches sign as $s$ is zero at a Legendrian point.  Consider now
a circle passing through such a Legendrian point. Then (see figure \ref{Fig:switchsign} (b)) if the
circle is to lie outside the occluded region, then moving in the
direction of the parametrization of $\gamma$, we must have that $s$
moves from positive values to negative values.  For $s$ to change
sign, we must have that the inward normal to the rule $c$ changes
direction.  As we have seen, this can happen at a Legendrian point,
but it may also happen if $|R|$ tends to either zero or infinity.  By
\rfL[CL]{Rnot0}, $R$ cannot be zero so we must have that either we
encounter another Legendrian point of this type or $|R| \ra \infty$.
Without loss of generality, we may pick the two Legendrian points so
that one of the two arcs of $\gamma$ that they bound contains no other
Legendrian points of this type.  Denote this region by $\gamma|_{(t_1,t_2)}$.  Hence, we must have that $|R| \ra
\infty$ for some $t_0\in(t_1,t_2)$.  

To finish the proof, suppose there exists a point $x$ in the occluded
region so that the projection of $\Sigma$ misses $x$.  For each $t$,
there exists a circle, $\overline{c}_t(s)$, connecting $\gamma(t)$ to
$x$ with $\overline{c}_t(0)=\gamma(t)$ and $\overline{c}'_t(0)=\frac{\nabla I}{|\nabla I|}^\perp(\gamma(t))$.  Let
$\overline{R}(t)$ be the radius of this circle and 
\begin{equation*}
\mathscr{R}(t) =
\begin{cases}
\overline{R}(t) \; \; \text{ if $\dot{\gamma}(t) \cdot
    \ddot{\overline{c}}_t(0)>0$}\\
-\overline{R}(t) \; \; \text{ if $\dot{\gamma}(t) \cdot
    \ddot{\overline{c}}_t(0)<0$}\\
\end{cases}
\end{equation*}

So, for the projection of $\Sigma$ to miss $x$, we must have that
$\text{sign}(\cos(\alpha))R(t) \neq \mathscr{R}(t)$ for all $t\in
  (t_1,t_2)$. But, by the discussion above, $
  \text{sign}(\cos(\alpha))R(t)$ changes sign on this region and
    hence $\overline{R}(t)$ must also tend to $\infty$ and so must tend to
    $\infty$ at $t_0$ as well.  By construction, the rule through
    $\gamma(t_0)$ is a straight line and hits $x$, contradicting the
    assumption that the circles at the Legendrian points lie outside $D$ (see figure \ref{Fig:switchsign} (c)).
\epf

\bgL{NoGaps}
If $Q^\prime <0$ everywhere then the limiting rules at any Legendrian point are external to the occluded disc.
\enL

\begin{figure}
\begin{center}
\includegraphics[width=2in]{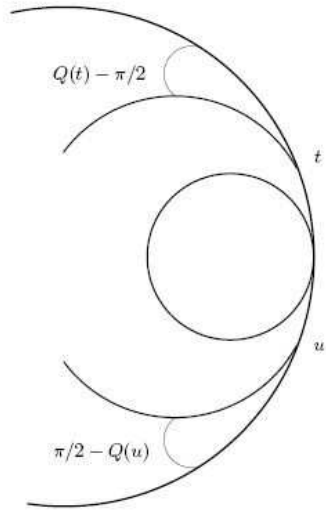}
\caption{Legendre Gap.}
\label{Fig:NoGap}
\end{center}
\end{figure}

\pf
Suppose not. Then without loss of generality we may rotate and reflect
the image data to match Figure \ref{Fig:NoGap}, where we are assuming
that the curves are oriented to the counter-clockwise
direction. Elementary arguments then show that the angles marked are
indeed $Q(t)-\pi/2$ and  $\pi/2-Q(u)$ and that these must therefore
both be positive. However at the Legendrian point we must have
$Q=\pi/2$. This clearly violates the condition that $Q^{\prime}<0$. 

\epf

\bgR{rem3}  We remark that if $Q'$ is positive at a Legendrian point, then the circle tangent to the occlusion boundary may indeed lie inside the occluded circle.
\enR

\bgT{goodcaseA}  Let $I: \R^2 \ra \R$ be an intesity function of an
image with an occlusion given by a circular region $D$.  Further,
suppose $\gamma \in \mathcal{RT}$ is the $\theta$ lift of $\partial
D$ and that the occlusion is completely nondegenerate and occludes
no critical points of $I$.  If $Q'(t) \neq 0$ for $t \in [0,2\pi]$
then there exists a minimal spanning surface of $\gamma$ satisfying
condition \ref{cond0}.  Moreover, if $Q'(t) <0$ for $t \in [0,2\pi]$ then the minimal surface satifies condition \ref{condA}.  
\enT

\pf As discussed above, if the occlusion is completely nondegenerate
and occludes no critical points then $deg \; Q =-1$.  Further, since
$Q' \neq 0$, we have that the function $u(t)$, implicitly defined by 
\[Q(t)+Q(u) = (2k+1)\pi\]
is monotone.  Thus, $u:S^1 \ra S^1$ is one to one and onto and thus,
for each $t$, there is a unique point connecting to to $t$ given by
$u(t)$.  Using the matrix equation \ref{mateqn} determines the rules
joining $t$ to $u(t)$ and, by picking the portion of the circle
defining the rule to be inside of $D$ when projected to $\R^2$, we
also satisfy condition \ref{cond0}.  

By \rfL{surj1} and \rfL{NoGaps}, we have that under the assumption of $Q'<0$, the
spanning surface satisfies condition \ref{condA}.

\epf

\section{Discussion} 
\begin{figure}
\setcounter{subfigure}{0}
\renewcommand{\thesubfigure}{1.\alph{subfigure} \space}
\fbox{\subfigure[Image]{\includegraphics[width=1.5in]{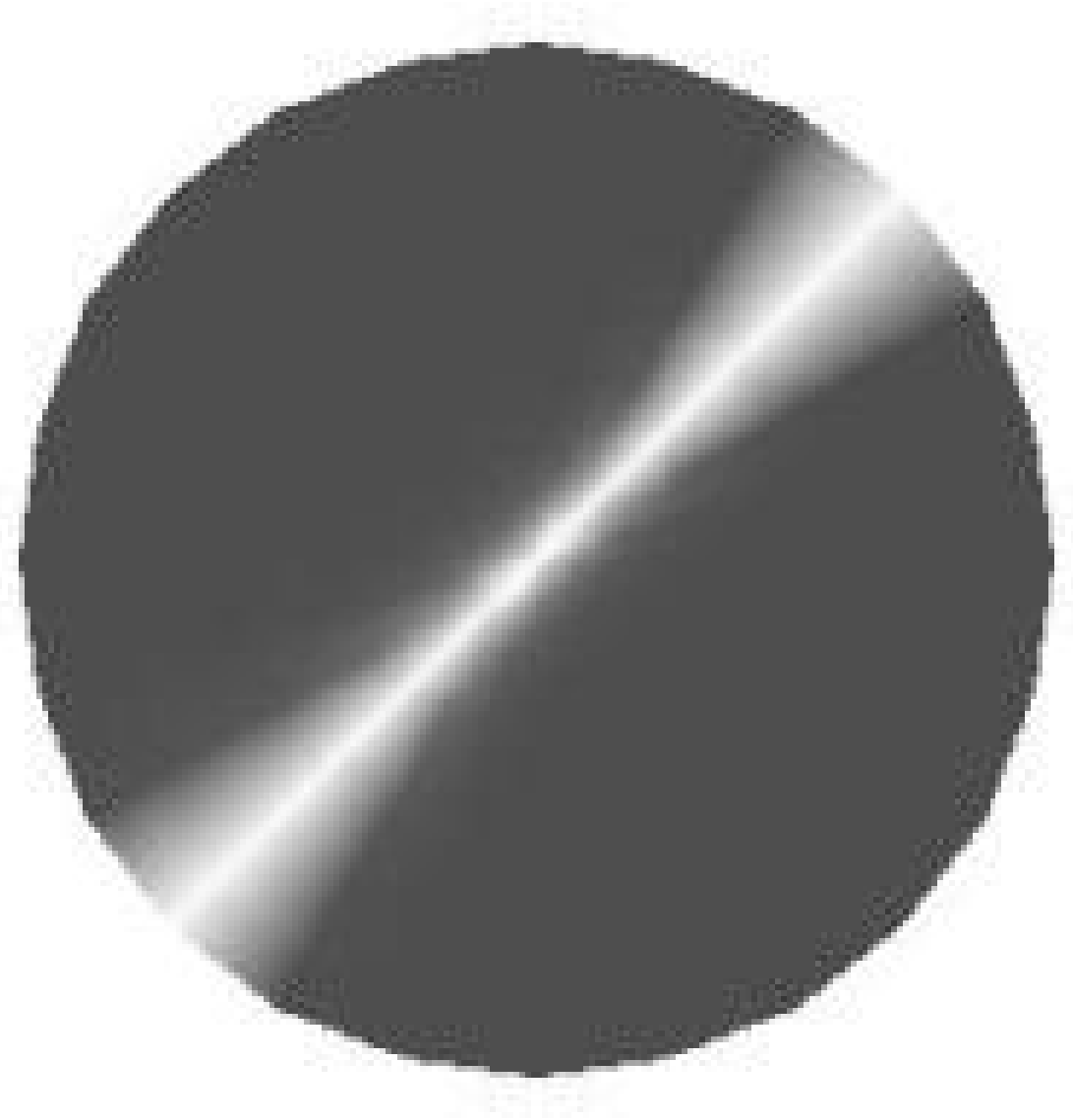}}
  \quad
  \subfigure[Occlusion]{\includegraphics[width=1.5in]{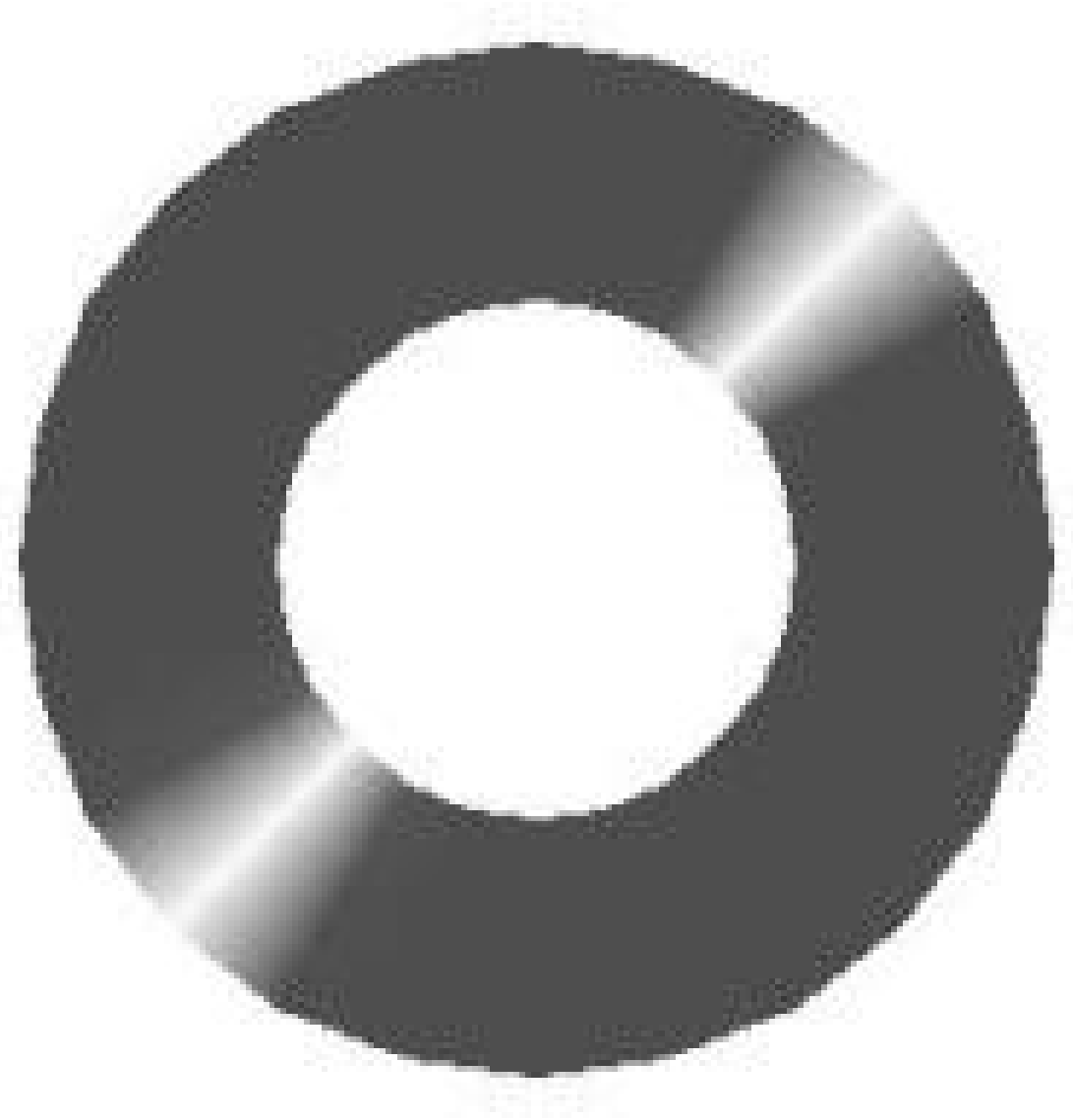}}
  \quad
  \subfigure[Completion]{\includegraphics[width=1.5in]{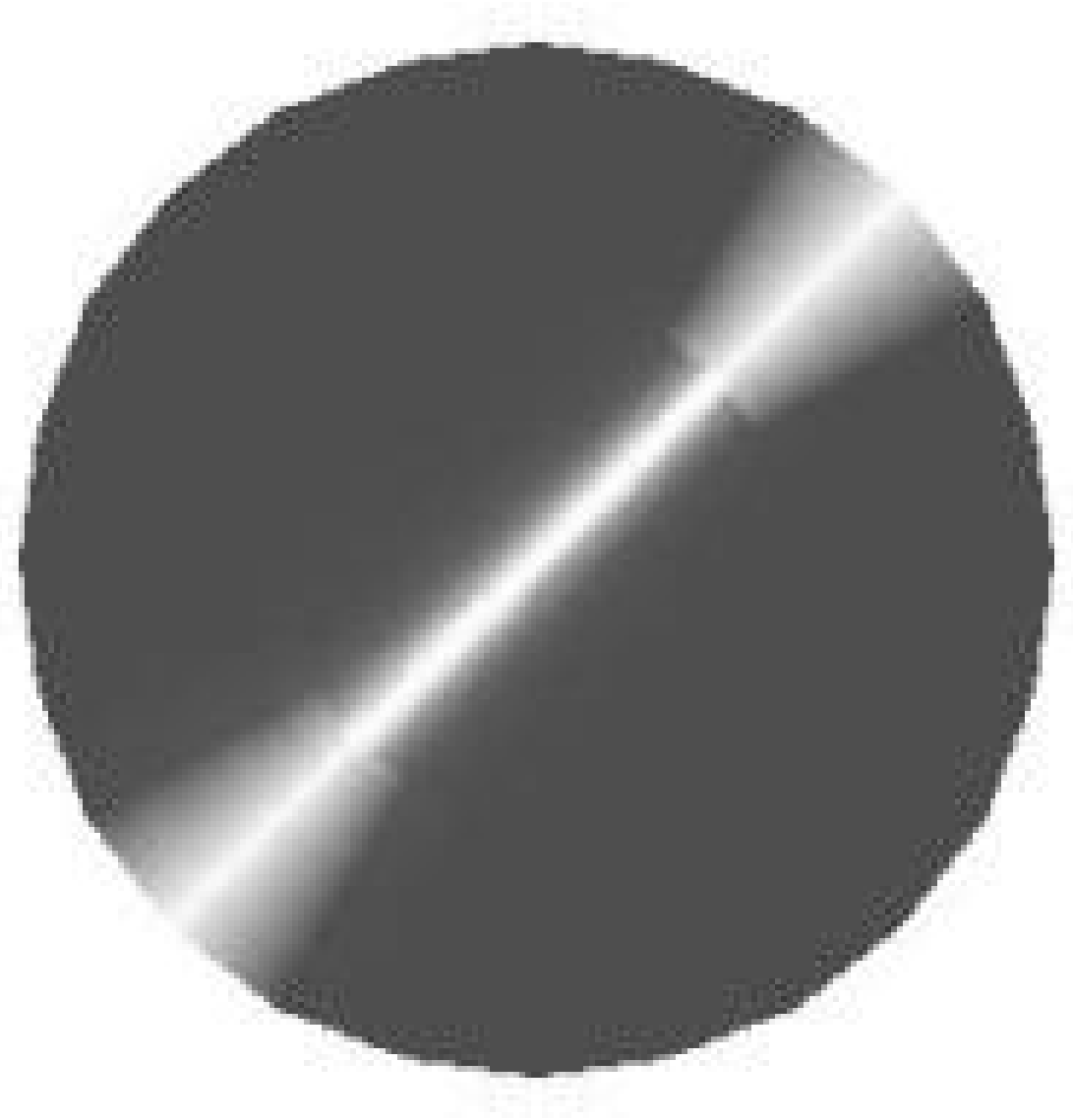}}}
\caption{Occlusion of a simple linear image}\label{occl1}
\end{figure}

\begin{figure}
\setcounter{subfigure}{0}
\renewcommand{\thesubfigure}{2.\alph{subfigure} \space}
\fbox{\subfigure[Image]{\includegraphics[width=1.5in]{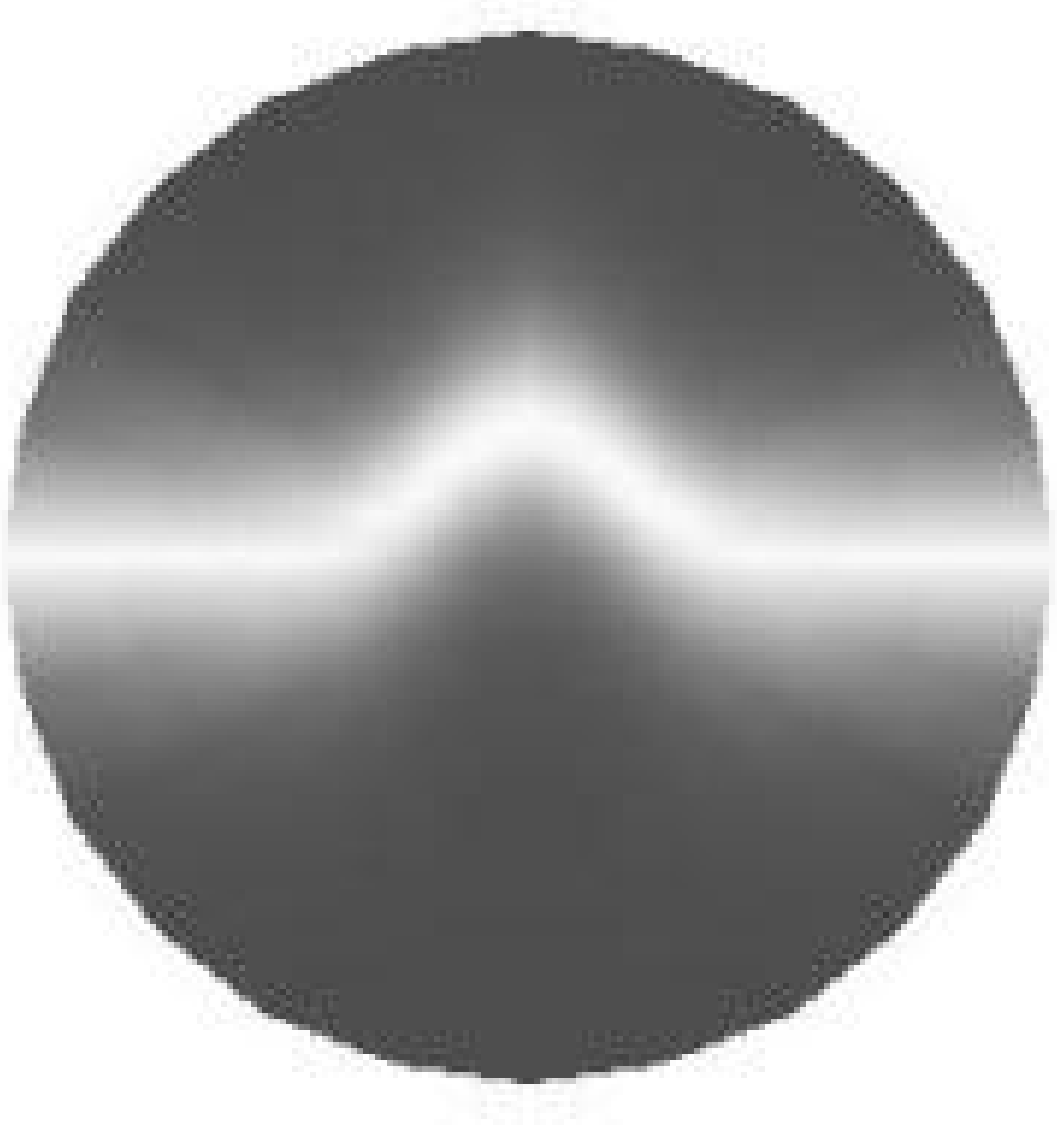}} \quad \subfigure[Occlusion]{\includegraphics[width=1.5in]{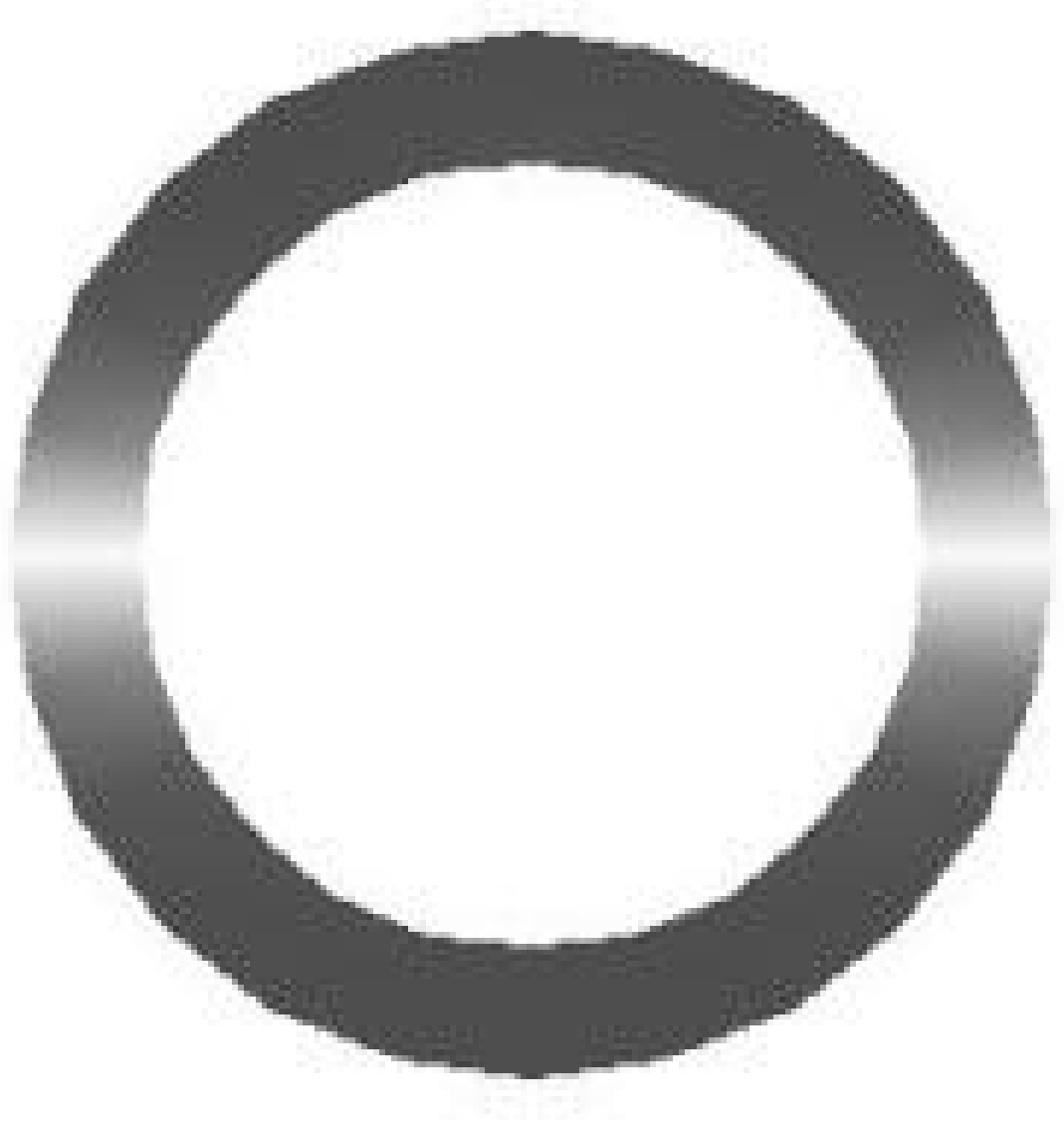}} \quad \subfigure[Completion]{\includegraphics[width=1.5in]{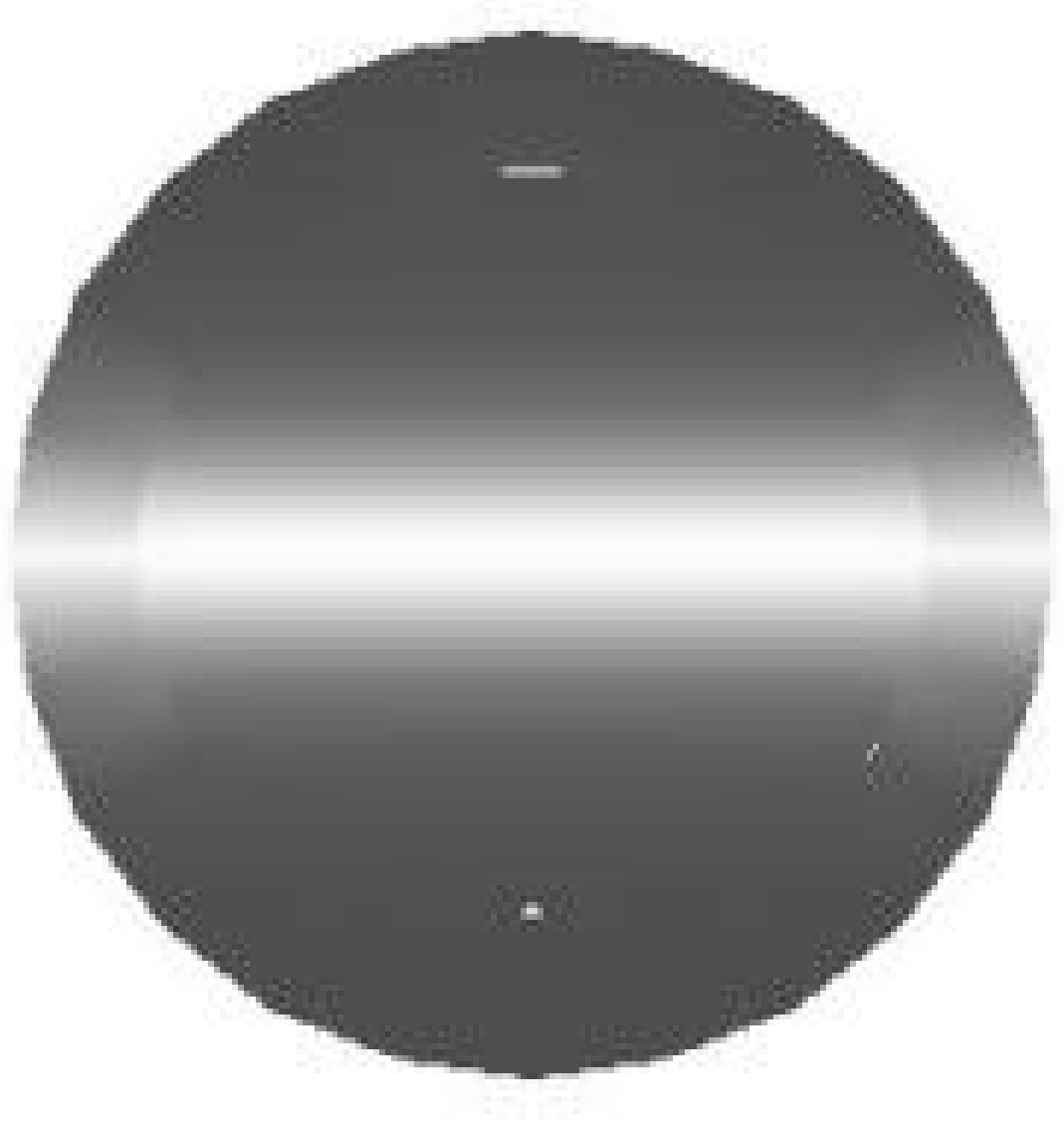}}}
\caption{Occlusion of a complicated portion of an image.}\label{occl2}
\end{figure}

\begin{figure}
\setcounter{subfigure}{0}
\renewcommand{\thesubfigure}{3.\alph{subfigure} \space}
\fbox{\subfigure[Image]{\includegraphics[width=1.5in]{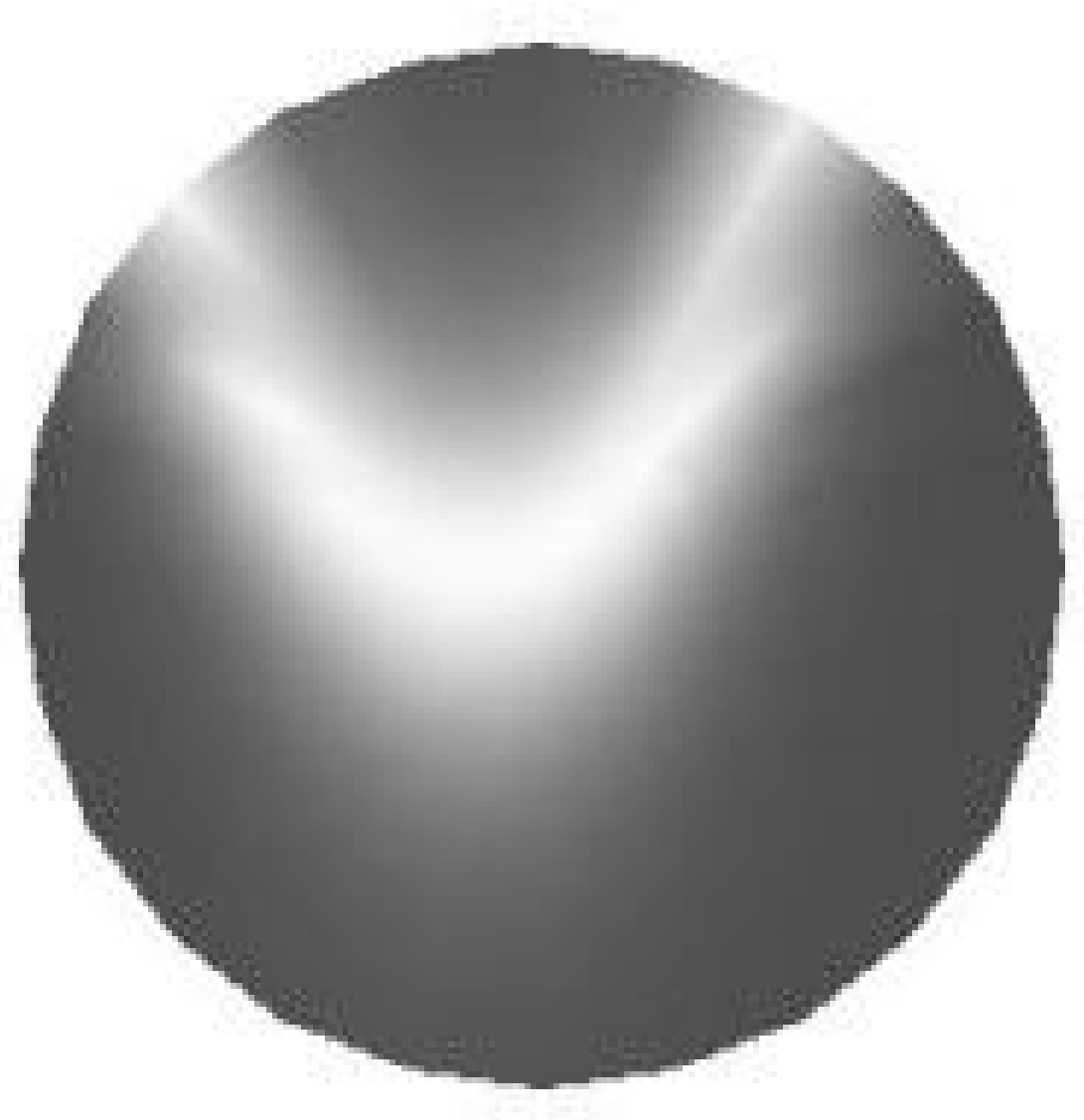}}
  \quad
  \subfigure[Occlusion]{\includegraphics[width=1.5in]{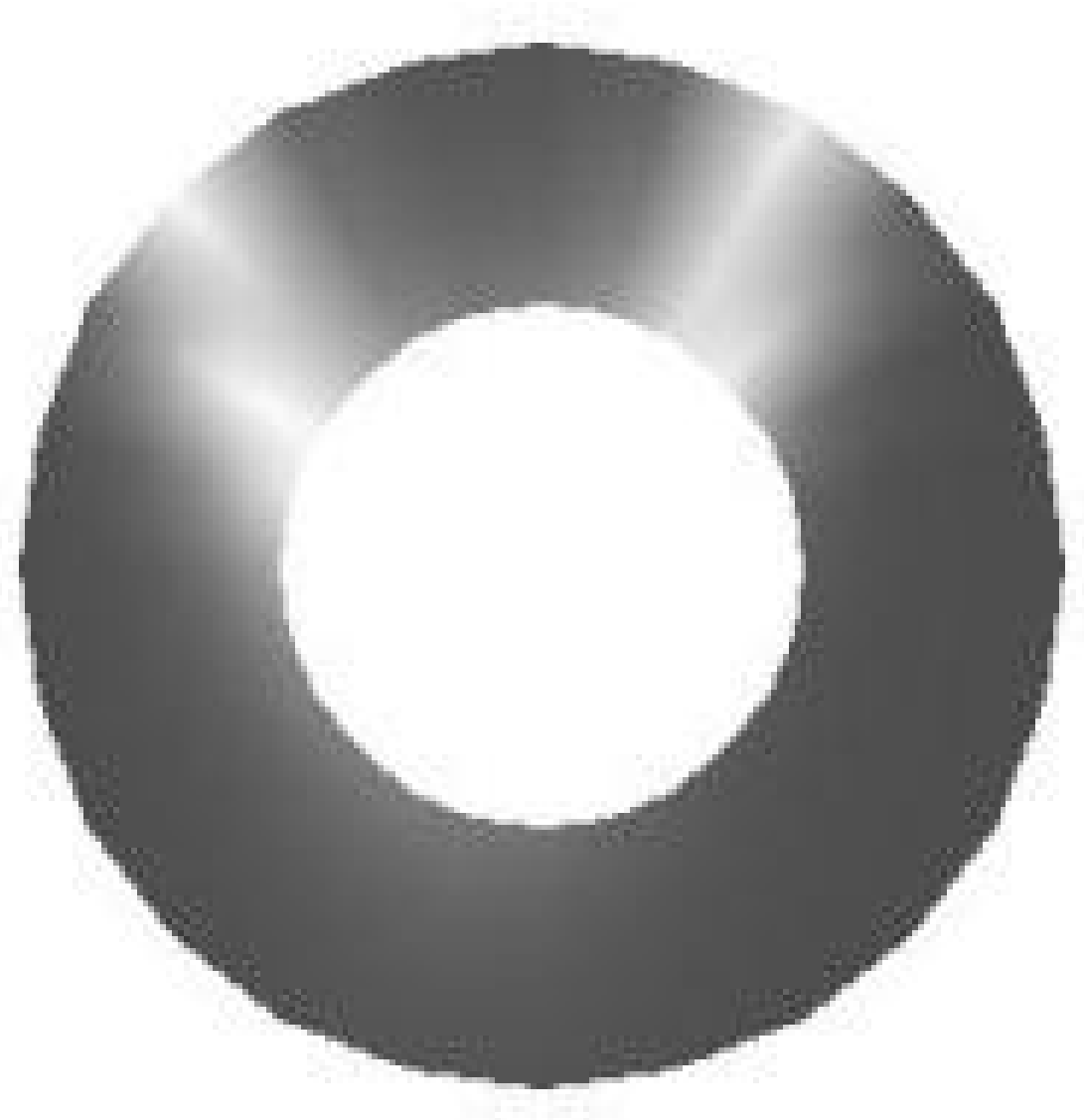}}
  \quad
  \subfigure[Completion]{\includegraphics[width=1.5in]{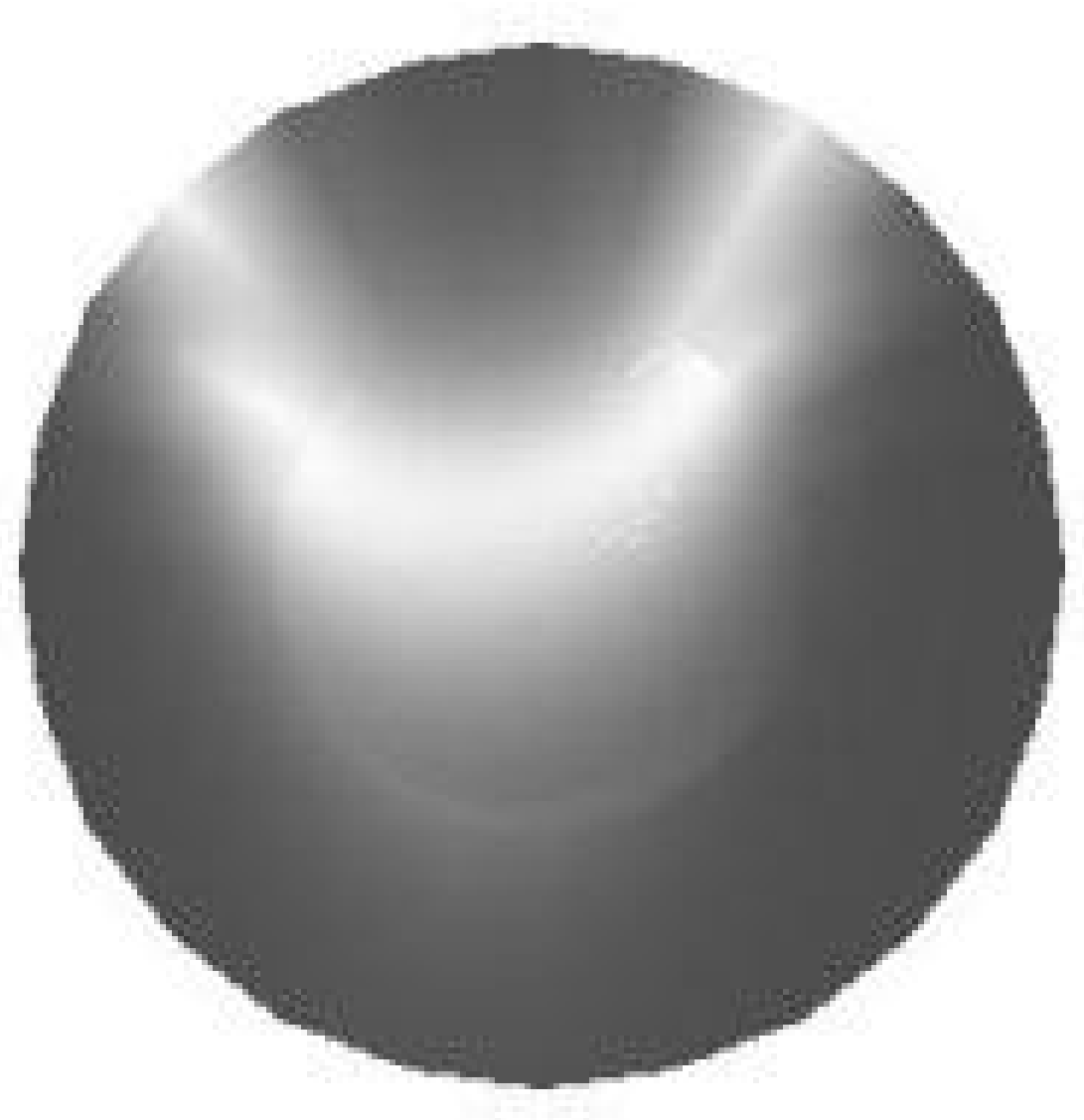}}}
\caption{Occlusion of a curve}\label{occl3}
\end{figure}

\begin{figure}
\setcounter{subfigure}{0}
\renewcommand{\thesubfigure}{4.\alph{subfigure} \space}
\fbox{\subfigure[Image]{\includegraphics[width=1.5in]{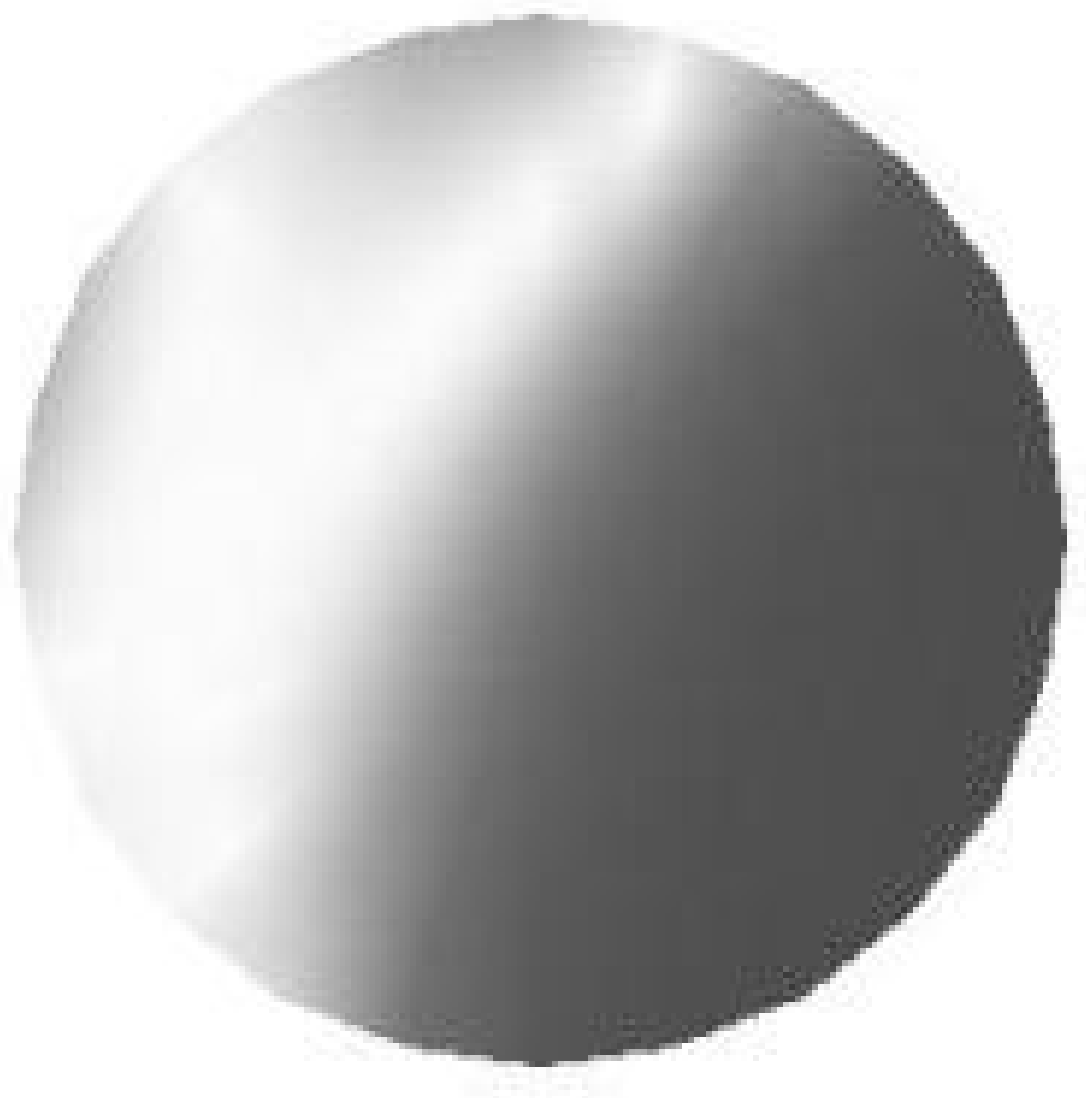}} \quad \subfigure[Occlusion]{\includegraphics[width=1.5in]{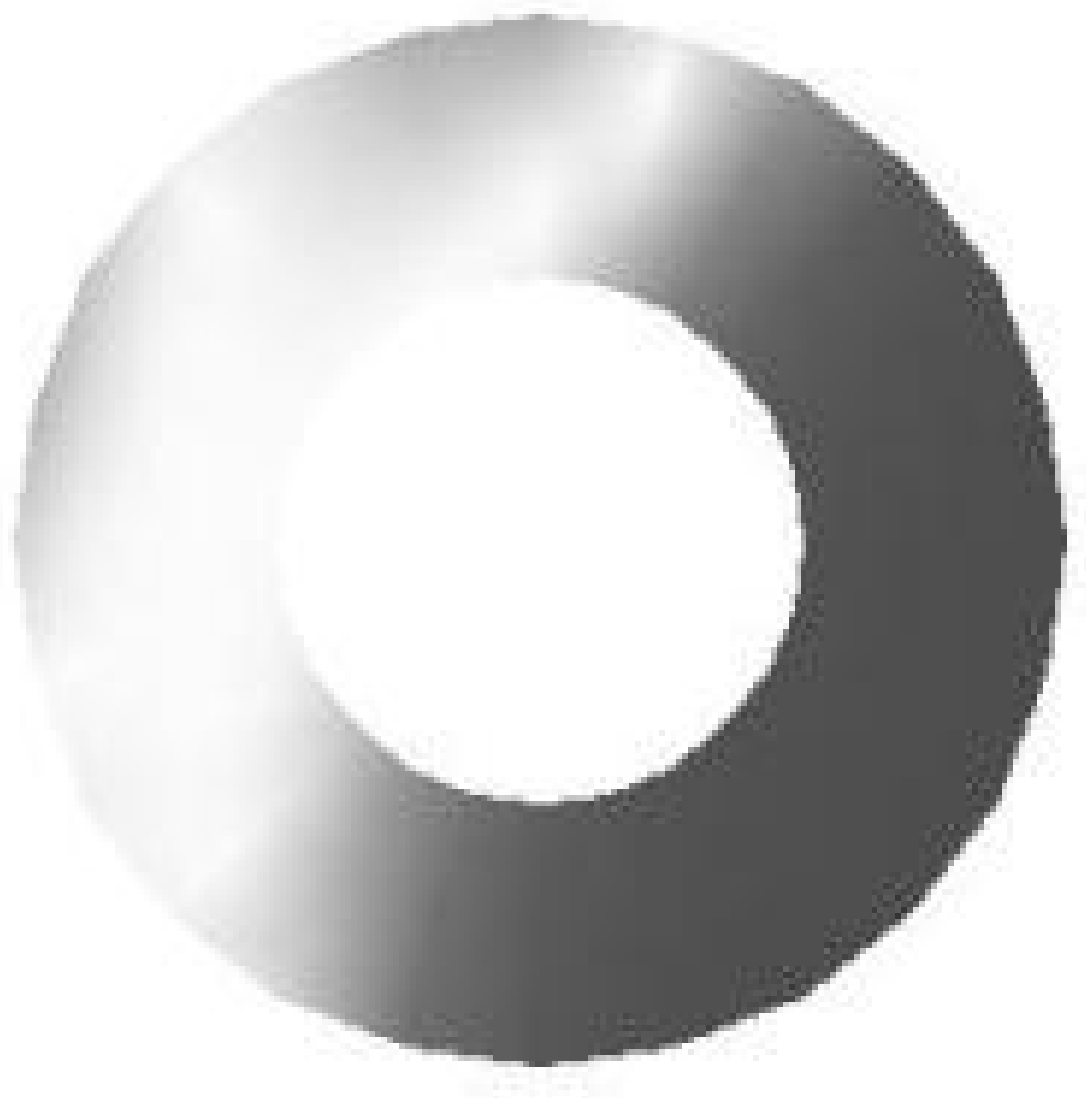}} \quad \subfigure[Completion]{\includegraphics[width=1.5in]{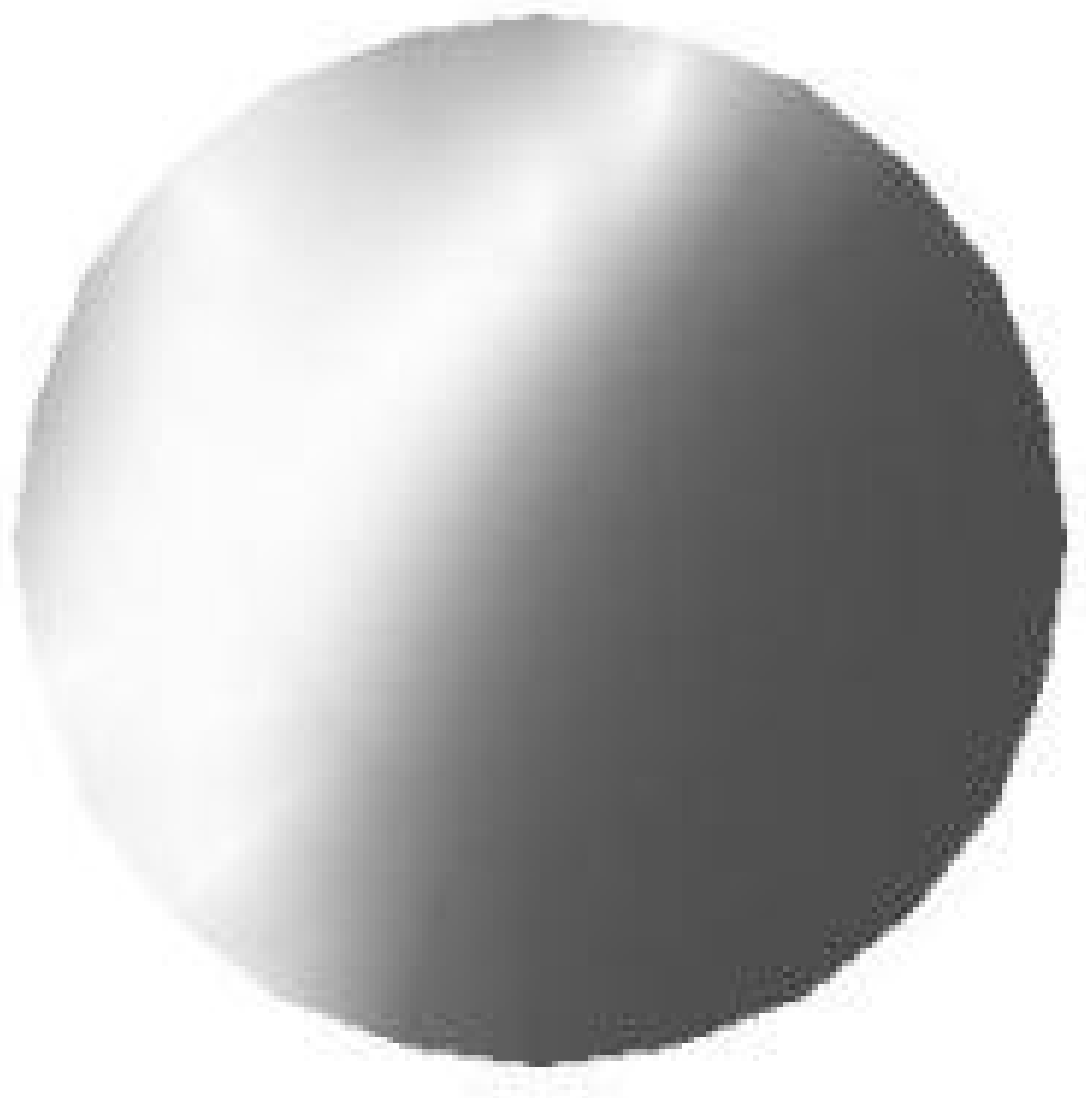}}}
\caption{Completion preserves concavity, I}\label{occl4}
\end{figure}

\begin{figure}
\setcounter{subfigure}{0}
\renewcommand{\thesubfigure}{5.\alph{subfigure} \space}
\fbox{\subfigure[Image]{\includegraphics[width=1.5in]{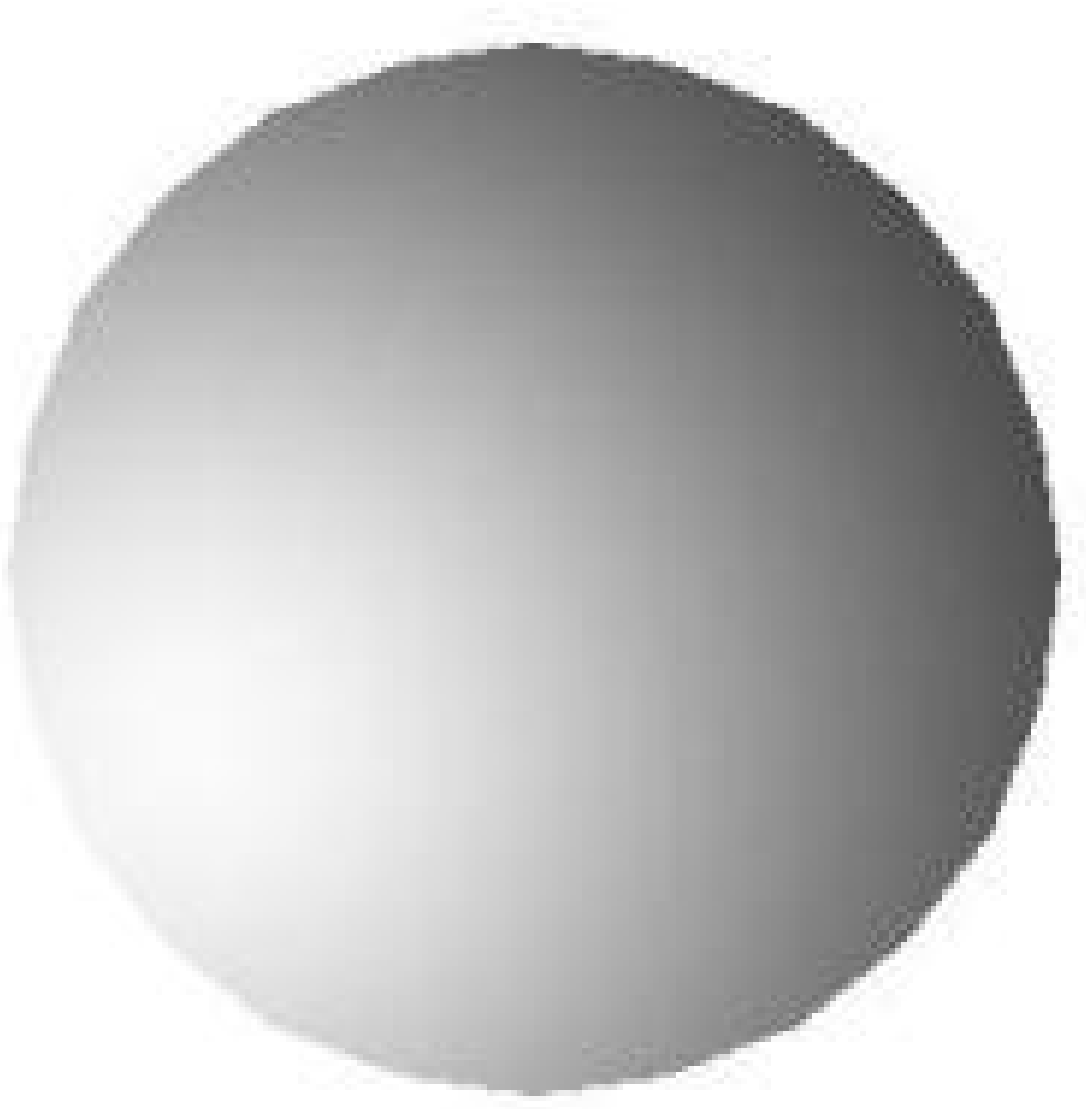}}
  \quad
  \subfigure[Occlusion]{\includegraphics[width=1.5in]{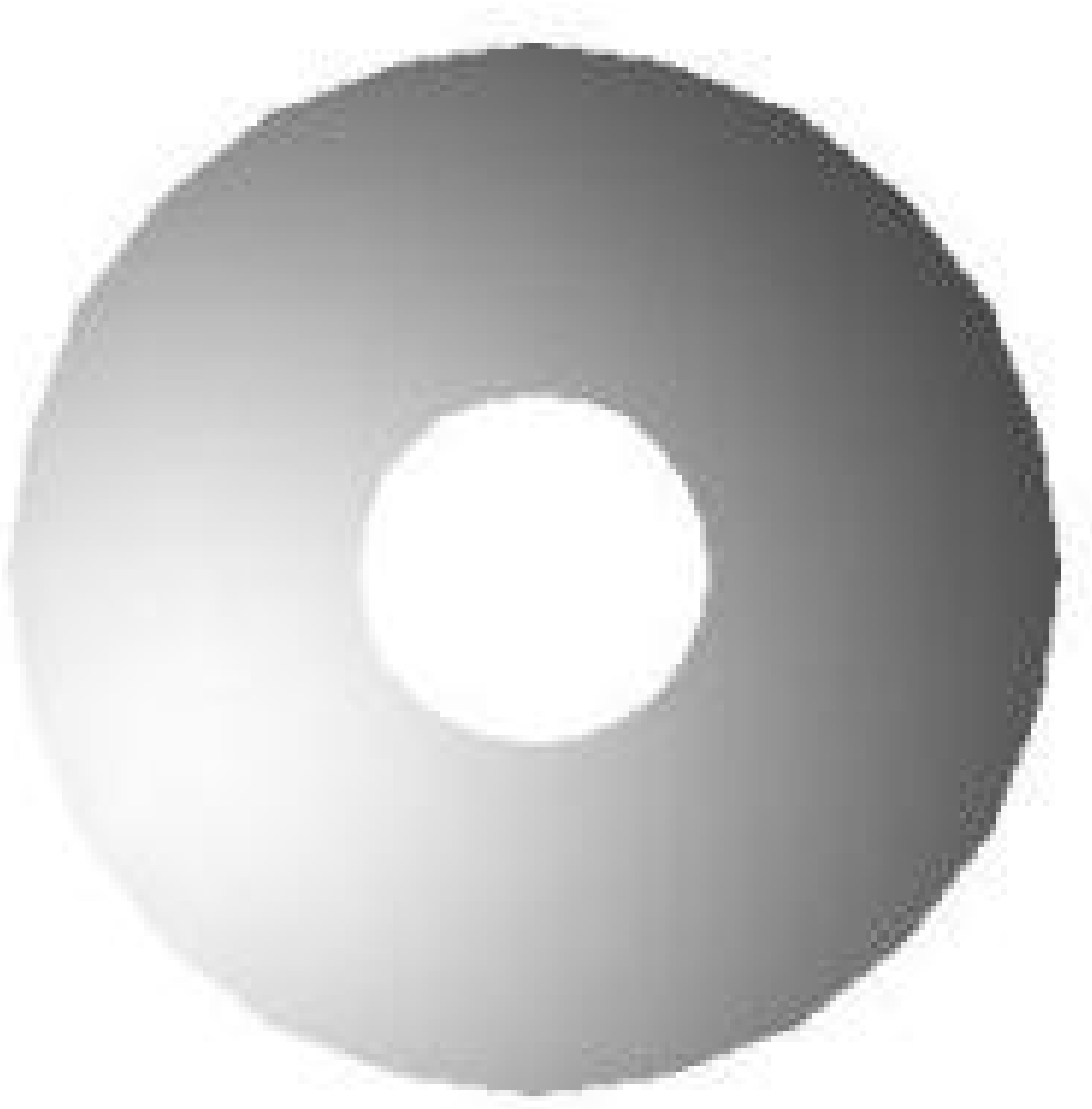}}
  \quad
  \subfigure[Completion]{\includegraphics[width=1.5in]{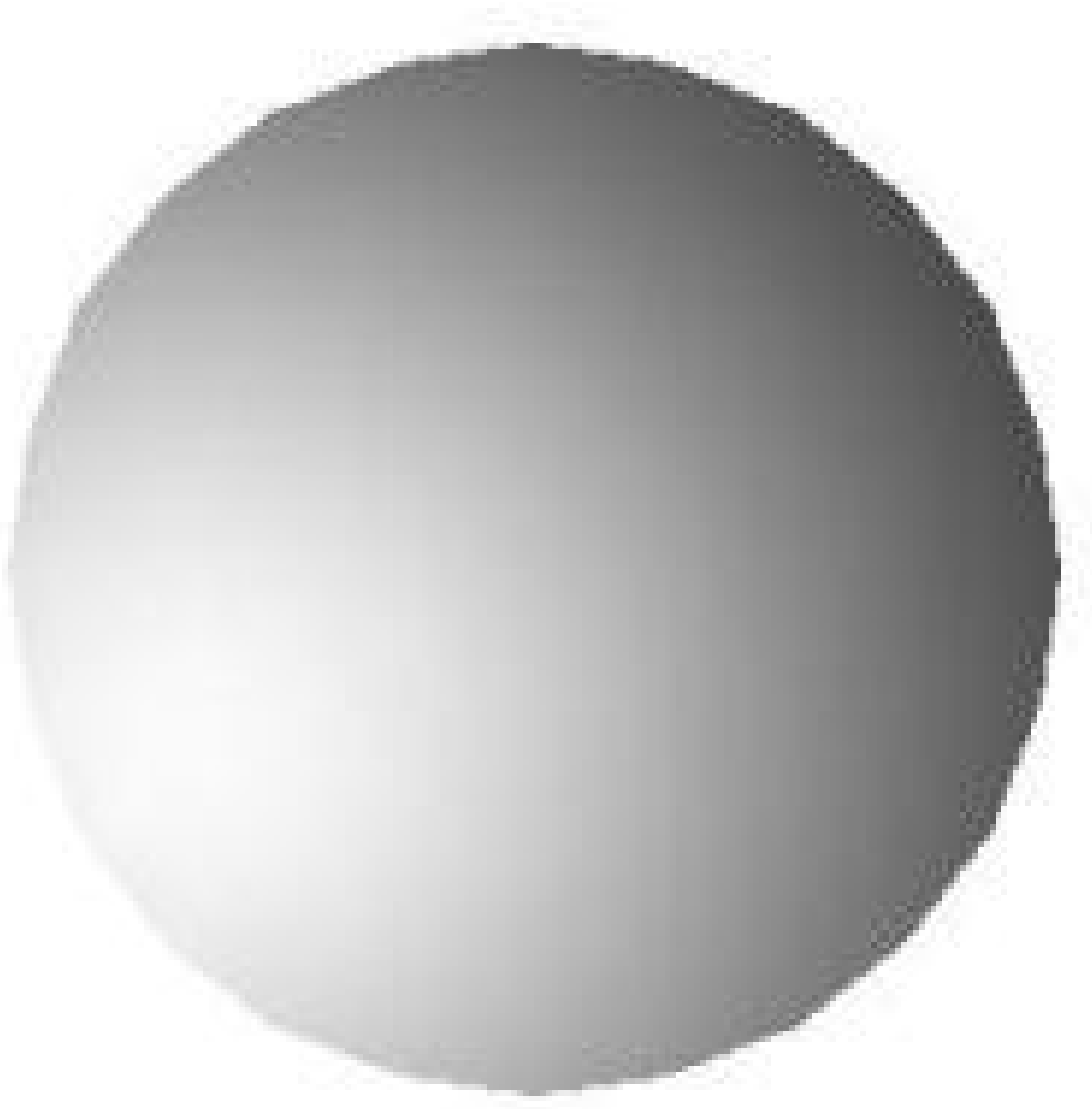}}}
\caption{Occlusion preserves concavity, II}\label{occl5}
\end{figure}

\setS{discuss}
The theorem of the previous section shows that if the behavior of $I$
near the occluded region is relatively tame, then we can easily
construct a minimal surface that spans the occluded region.  In
figures \ref{occl1} - \ref{occl5}, we demonstrate an implementation of
the algorithm described in the previous section in cases where the
theorem applies.  The images were created by assigning color to the
value of various functions (as with the previous test image).  


In figure \ref{occl1}, we see that the algorithm easily completes a
linear image using straigt lines.  Similarly, in figure \ref{occl2},
we see that when the a complicated portion of a region is occluded, it
may not necessarily be recovered.  In more interesting examples that
show the power of this new method, figures \ref{occl3} ,\ref{occl4},
and \ref{occl5}
show completion using the circles idenitified in previous section.
Figure \ref{occl3} shows the completion of a curve which preserves the
curvature of the original image.  Figures \ref{occl4} and \ref{occl5} show that the
algorithm preserves concavity.  We note that this are features that
earlier diffusion based algorithms often had trouble completing in a
reliable manner.

All of the other cases presented show different types of pathology:
\begin{enumerate}
\item Case 2 shows the simplest type of failure of a spanning surface
  to be a graph:  when $Q'$ has zeros.  In this case, we have
  ``backtracking'' of rules which causes the resulting surface derived
  from the implicitly defined function $u$ to give an immersed rather
  than embedded surface.  
\item Case 3 shows an instance where no smooth minimal surface
  exists due to the presence of non-Legendrian solitary points.  As
  pointed out above, this behavior seems to be characteristic of
  occluded critical points of $I$.  Moreover, this behavior further
  indicates that restricting to smooth spanning surfaces, while
  computationally effective, will not solve any possible minimal
  surface problem with Dirichlet conditions.  Case 2 is  similar to a result
  to the case in the Heisenberg group\cite{Pauls:Obstr} but Case 3
  shows entirely new behavior stemming from the nontrivial topology of
  the space.  
\item Case 4 shows that with higher degree there are potentially both
  nonuniqueness issues as well as problems with satisfying condition
  \ref{condB}, i.e. forcing the spanning surface to be a graph.
\end{enumerate}

The pathologies outlined above, coupled with the discussion of connection $\gamma$ to
  $\overline{\gamma}$, lead us to several conclusions.  First, the
  restriction to smooth spanning surfaces, while sufficient for many
  types of problems (such as those of theorem \rfT[thms]{goodcaseA}) is likely
  insufficient for more complicated areas of an image.  Second, the
  different types of pathologies suggest that to amodally complete a
  given occlusion, the ``best'' completion is likely to come from
  knitting together various pieces of several different solutions
  (i.e. from different branches of the curve defining $u$ or from
  pieces connecting $\gamma$ to $\overline{\gamma}$).  Again, this
  points towards the neccesity of a more sophisticated mechanism.
  However, we point out that this is consistent with the simulation
  data found by Citti and Sarti\cite{CS} showing that several
  different possible completions are present at the same time in
  $\mathcal{RT}$ after using their diffusion method.  If this model of
  minimal surface completion is accurate reflecting the completion
  mechanism in V1, this ambiguity stemming from multiple (partial)
  solutions may be resolved by the input and feedback from other
  layers fo the visual cortex.  In particular, we note that our
  algorithm often produces connections between level sets of different
  ``heights'' thus creating a completion which is not ideal from the
  point of view of matching like intensities within the image.  This
  is consistent with the model of the visual cortex present in section
  \ref{model} as the representation of the image in $\mathcal{RT}$ does
  not carry information about the intensity of the image, but only
  information about the level sets of the image itself.  One expects
  that with additional input such as color/intensity information, the
  best possible completion could be picked out of the possibilities.

\section*{Acknowledgment}
Both authors are partially supported by NSF grant DMS-0306752.


\begin{thebibliography}{10}

\bibitem{V:AM}
L.~Ambrosio and S.~Masnou.
\newblock A direct variational approach to a problem arising in image
  reconstruction.
\newblock {\em Interfaces and Free Boundaries}.
\newblock To appear.

\bibitem{CH}
Jih-Hsin Cheng and Jenn-Fang Hwang.
\newblock Properly embedded and immersed minimal surfaces in the {H}eisenberg
  group.
\newblock 2004.
\newblock Preprint: arxiv math.DG/0407094.

\bibitem{CHMY}
Jih-Hsin Cheng, Jenn-Fang Hwang, Andrea Malchiodi, and Paul Yang.
\newblock Minimal surfaces in pseudohermitian geometry.
\newblock 2003.
\newblock Preprint.

\bibitem{CMS}
G.~Citti, M.~Manfredini, and A.~Sarti.
\newblock Neuronal oscillation in the visual cortex: {G}amma-convergence to the
  {R}iemannian {M}umford-{S}hah functional.
\newblock {\em SIAM Jornal of Mathematical Analysis}, 35(6):1394 -- 1419, 2003.

\bibitem{CS}
G.~Citti and A.~Sarti.
\newblock A cortical based model of perceptual completion in the
  roto-translation space.
\newblock 2004.
\newblock Preprint.

\bibitem{Cole}
Daniel Cole.
\newblock {\em On minimal surfaces in {M}artinet-type spaces.}
\newblock PhD thesis, Dartmouth College, 2005.

\bibitem{DGN}
D.~Danielli, N.~Garofalo, and D.-M. Nhieu.
\newblock {M}inimal surfaces, surfaces of constant mean curvature and
  isoperimetry in {C}arnot groups.
\newblock August, 2001.
\newblock Preprint.

\bibitem{GarNh}
Nicola Garofalo and Duy-Minh Nhieu.
\newblock Isoperimetric and {S}obolev inequalities for
  {C}arnot-{C}arath\'eodory spaces and the existence of minimal surfaces.
\newblock {\em Comm. Pure Appl. Math.}, 49(10):1081--1144, 1996.

\bibitem{GP}
Nicola Garofalo and Scott~D. Pauls.
\newblock The {B}ernstein problem in the {H}eisenberg group.
\newblock 2003.
\newblock Submitted.

\bibitem{V:GDIW}
C.D. Gilbert, A.~Das, M.~Ito, and G.~Westheimer.
\newblock Spatial integration and cortical dynamics.
\newblock {\em Proceedings of the National Academy of Sciences USA},
  93:615--622.

\bibitem{HP3}
Robert Hladky and Scott~D. Pauls.
\newblock A disocclusion algorithm based on a model of the visual cortex.
\newblock 2005.
\newblock In preparation.

\bibitem{Pauls:cmcgroup}
Robert~K. Hladky and Scott~D. Pauls.
\newblock Constant mean curvature surfaces in sub-riemannian spaces.
\newblock 2005.
\newblock Preprint.

\bibitem{V:Hoffman}
William~C. Hoffman.
\newblock The visual cortex is a contact bundle.
\newblock {\em Appl. Math. Comput.}, 32(2-3):137--167, 1989.
\newblock Mathematical biology.

\bibitem{V:HW62}
D.~H. Hubel and T.~N. Weisel.
\newblock Receptive fields, binocular interaction and functional architecture
  in the cat's visual cortex.
\newblock {\em J. Physiol.}, 160:106--154, 1962.

\bibitem{V:HW77}
D.~H. Hubel and T.~N. Weisel.
\newblock Functional architecture of macaque monkey visual cortex.
\newblock {\em Proc. R. Soc. London (Biol)}, 198:1--59, 1977.

\bibitem{Kobayashi}
S.~Kobayashi and K.~Nomizu.
\newblock {\em Foundations of Differential Geometry}.
\newblock John Wiley \& Sons, Inc., 1963.

\bibitem{LM}
G.~P. Leonardi and S.~Masnou.
\newblock On the isoperimetric problem in the {H}eisenberg group
  {${\mathbf{H^n}}$}.
\newblock Preprint, 2002.

\bibitem{LR}
G.~P. Leonardi and S.~Rigot.
\newblock Isoperimetric sets on {C}arnot groups.
\newblock {\em Houston J. Math.}, 29(3):609--637 (electronic), 2003.

\bibitem{V:NMS}
M.~Nitzberg, D.~Mumford, and T.~Shiota.
\newblock {\em Filtering, segmentation and depth}, volume 662 of {\em Lecture
  Notes in Computer Science}.
\newblock Springer-Verlag, Berlin, 1993.

\bibitem{Pauls:Obstr}
Scott~D. Pauls.
\newblock H-minimal graphs of low regularity in the {H}eisenberg group.
\newblock {\em Comm. Math. Helv.}
\newblock to appear.

\bibitem{Pauls:minimal}
Scott~D. Pauls.
\newblock Minimal surfaces in the {H}eisenberg group.
\newblock {\em Geom. Ded.}, 104:201--231, 2004.

\bibitem{V:P}
J.~Petitot.
\newblock The neurogeometry of pinwheels as a sub-{R}iemannian contact
  structure.
\newblock {\em J. Physiology}, 97:265--309, 2003.

\bibitem{V:PT}
J.~Petitot and Y.~Tondut.
\newblock Vers une neuro-geometrie. fibrations corticales, structures de
  contact et contours subjectifs modaux.
\newblock {\em Mathematiques, Informatique et Sciences Humaine, EHESS, Paris},
  145:5--101, 1998.

\bibitem{RR}
Manuel Ritor\'e and C\'esar Rosales.
\newblock Rotationally invariant hypersurfaces with constant mean curvature in
  the {H}eisenberg group $\mathbb{H}^n$.
\newblock 2005.
\newblock Preprint.

\bibitem{Tanaka}
N.~Tanaka.
\newblock {\em A differential geometric study on strongly pseudoconvex
  manifolds}.
\newblock Kinokuniya Book-Store Co., Ltd., 1975.

\bibitem{Webster}
S.M. Webster.
\newblock {Pseudo-Hermitian structures on a real hypersurface}.
\newblock {\em J. Differential Geometry}, 13:25--41, 1978.

\end{thebibliography}

\end{document}